\documentclass[12pt]{amsart}
\pdfoutput=1
\usepackage{geometry}                
\usepackage{graphicx}
\usepackage{amssymb}
\usepackage{epstopdf}
\usepackage[author-year]{amsrefs}  

\newtheorem{theorem}{Theorem}

\def\Z{\mathbb Z}

\newcommand{\comment}[1]{} 

\comment{
Files to include with this .tex file
DrawingButterflies.bbl  (the bibliography needed by arXiv. For me, I need the .bib file)
mathieu50_sym.pdf
Butterfly50_10.pdf
BF_t_1+.pdf
BF_t_1-.pdf
BF_t_2a_jump.pdf
BF_t_2b_jump.pdf
BF_t_2a.pdf
BF_t_2b.pdf
BF_t_3a_jump.pdf
BF_t_3a.pdf
BF_t_3b_jump.pdf
BF_t_3b.pdf
BF_t_3c_jump.pdf
BF_t_3c.pdf
gap_check_4_lbl.pdf
BF_t_4+.pdf
BF_t_5+.pdf
gaplines2a.pdf
gaplines2b.pdf
gaplines3a.pdf
gaplines3b.pdf
gaplines3c.pdf
gaplines4.pdf
mathieu100_27.pdf
BF_2_7a.pdf
Butterfly50_5_16.pdf
Butterfly50_5.pdf
Butterfly50_5_30.pdf
Butterfly50_5.pdf
Butterfly50_20.pdf
}

\title{Drawing butterflies from the almost Mathieu operator}
\author{Michael P. Lamoureux}
\address[M.~Lamoureux]{Dept.\ Mathematics and Statistics \\ University of Calgary \\ 2500 University Ave NW \\ Calgary AB T2N 1N4 \\ Canada}
\email{mikel@math.ucalgary.ca}
\date{May 7, 2010}                     

\begin{document}

\begin{abstract}
Plotting spectra of a range of almost Mathieu operators reveals a beautiful fractal-like image that contains multiple copies of a butterfly image. We demonstrate that plotting the butterflies using a gap-labelling scheme based on K-theory or Chern numbers reveals systematic discontinuities in the gap positioning. A proper image is produced only when we take into account these discontinuities, and close the butterfly wingtips at the points of discontinuity. A conjecture is presented showing a simple formula for locating the discontinuities, and numerical evidence is given to support the conjecture. We also present new renderings of this butterfly.
\end{abstract}

\maketitle

\section{Introduction}

The inspiration for this work is the diagram in Figure~1, a fractal-like structure known as the Hofstadter butterfly~\cite{hof76}, which represents the energy levels of an electron travelling through a periodic lattice under the influence of a magnetic field; that is, a Bloch electron~\cite{brown64}. Each horizontal line in the image corresponds to union of spectra of certain bounded self-adjoint operators $H_{\theta,\phi}$ on the Hilbert space $l^2(\Z)$, defined by the equation
\[ (H_{\theta,\phi} \mathbf{x})_n = \mathbf{x}_{n-1} + \mathbf{x}_{n+1}
+ 2\cos(2\pi n \theta + \phi) \mathbf{x}_{n}, \]
where the parameter $\theta\in[0,1]$ corresponds to the vertical position in the image in Figure~1, and $\phi$ is a physical phase parameter.
This operator, known as the almost Mathieu operator  and also  as Harper's operator~\cite{last94}, is a discrete, one-dimensional Schr\"{o}dinger operator, which has been the object of study for many important physical problems~\cite{gold09}.
\begin{figure}[ht] 
   \centering
   \includegraphics[width=6in]{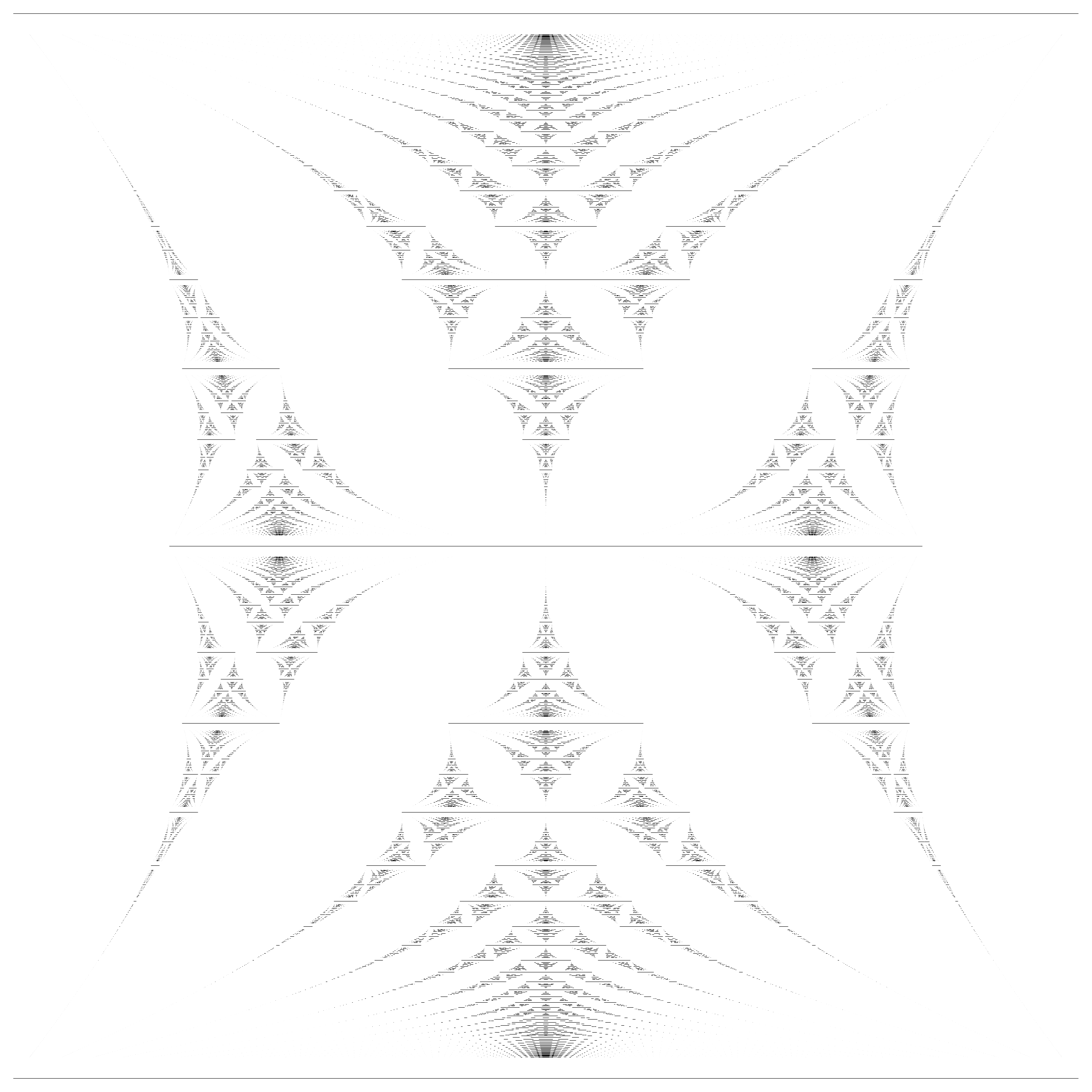} 
   \caption{The Hofstadter butterfly.}
   \label{fig:Hof}
\end{figure}

\begin{figure}[ht] 
   \centering
    \includegraphics[width=6in]{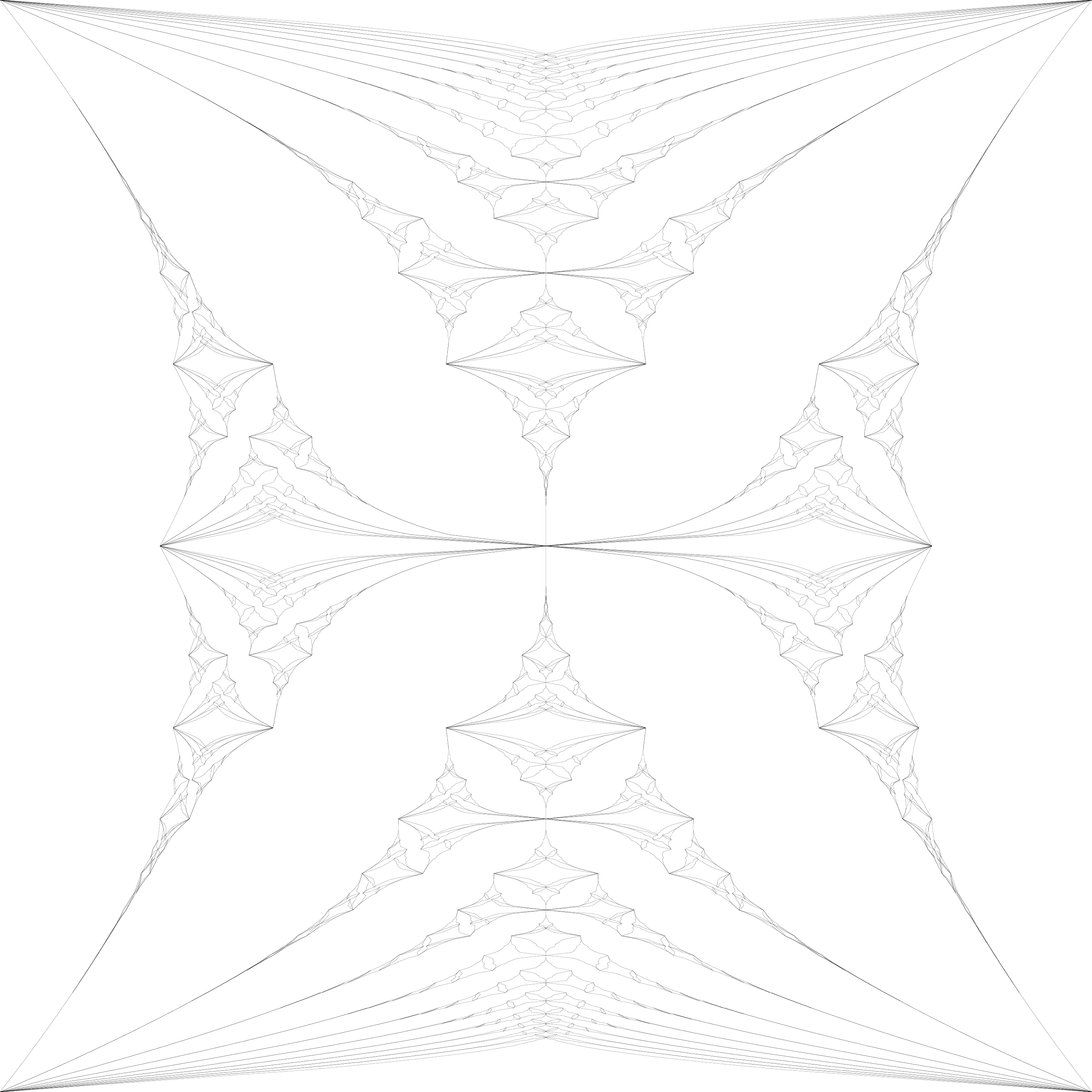} 
   \caption{Drawing only the butterflies.}
   \label{fig:BF50}
\end{figure}

Equivalently, each horizontal line in Figure~1 represents the spectrum of the self-adjoint operator 
\[ h_\theta = u + u^* + v + v^* \]
where unitaries $u,v$ are generators of the rotation C*-algebra $A_\theta$, characterized by a commutation rule
\[ vu = e^{2\pi i \theta} uv. \]
These particular algebras are known as non-commutative tori and are of considerable interest in noncommutative geometry~\cite{connes94}. The Cantor-like structure of the spectra for operators in these algebras has been a topic of study in C*-algebras for quite some time, as discussed in~\cites{avila06, bel82, choi90, last94, puig03}

Our interest in this paper is in the mechanics of drawing the image in Figure~1, its related image in Figure~2, and what this reveals about gap labelling theorems~\cites{bel90,gold09,kam03,ypma07} for these Schr\"{o}dinger operators. This particular image in Figure~1 was created in high resolution, using the numerical package MATLAB~\cite{matlab} to accurately compute the eigenvalues of certain tridiagonal matrices, then plotting line segments between eigenvalues using PostScript commands to ensure a high quality image. It is useful to refer to \cite{cass05} for details on how to create such high resolution mathematical illustrations in PostScript. 

In the finely rendered image in Figure~1, the viewer's eye is immediately drawn to butterfly-like structures in the diagram which repeat over and over at various scales and positions. Following Casselman's dictum that {\em mathematical drawings should represent concepts, not accidents}~\cite{cass05}, it is fascinating to consider how we might draw only the butterflies. In Figure~2, we do exactly that, drawing precisely the butterflies that are inherent in the image in Figure~1, without drawing all the spectral lines themselves. 

This new high resolution image was also created using Postscript commands, based on numerical calculations of spectra of operators generated in MATLAB. This new image highlights the fractal symmetry of the diagram and suggests further investigation on the explicit symmetries.

The goal of this paper is to record the algorithmic process used to image only the butterflies, using the gap labelling procedure for the almost Mathieu operator. We discover through this process that, from a certain point of view,  there are discontinuities in gap labelling, and we present a conjecture for specifying where these discontinuities occur. It appears this conjecture is new to the literature, but we are at a loss to prove it. The conjecture is verified for low rational values of parameter $\theta$.

A secondary goal is to reveal this new approach to rendering a high resolution image of the Hofstadter butterfly.

\section{Spectral calculations}

The spectral calculations for $h_\theta$ are well-known, and are particularly straight-forward when $\theta$ is rational~\cites{choi90,lam07}. For $\theta=p/q$ in reduced form, irreducible representations of the rotation C*-algebra $A_\theta$ are obtained using $q\times q$ matrices $U,V$, with $U$ a cyclic permutation matrix and $V$ a diagonal matrix with entries given as consecutive integer powers of $e^{2\pi i p/q}$.   Under the irreducible representation, the operator $h_\theta$ maps to the $q\times q$ matrix
\[ H_{\theta,z_1,z_2} = z_1U + \overline{z_1}U^* + z_2V + \overline{z_2}V^*, \]
where $z_1,z_2$ are complex numbers of modulus one, specifying the different irreducible representations. The spectrum of $h_\theta$ is just the union of the point spectra for the various matrices $H_{\theta,z_1,z_2}$ over the range of $z_1,z_2$

The characteristic polynomial of the $q\times q$ matrix $H_{\theta,z_1,z_2}$ is invariant under changes of parameters $z_1,z_2$, except for a zero-order term $z_1^q + z_1^{-q} + z_2^q + z_2^{-q}$. Consequently,  the union of spectra over the set of parameters $z_1,z_2$ is just the inverse image, under a degree q polynomial, of the range of values for the term $z_1^q + z_1^{-q} + z_2^q + z_2^{-q}$. This range is the interval $[-4,4]$, and the inverse image is a set of $q$ distinct intervals on the real line.

These $q$ intervals are all disjoint, except in the case $q$ even, in which case the central two intervals just touch at the spectral point zero. Thus, there are precisely $q-1$ gaps between intervals in the spectrum of $h_\theta$, where we include the ``empty gap'' at spectral value zero, for convenience.

To compute the spectral intervals, it suffices to compute the endpoints of each interval and then connect with a line. The endpoints of the interval are obtained by choosing values of $z_1,z_2$ such that the constant term $z_1^q + z_1^{-q} + z_2^q + z_2^{-q}$ obtains its extreme values $\pm 4$. 

The $q\times q$ matrix $H_{\theta,z_1,z_2}$ is almost tridiagonal -- and by choosing the extreme values for $z_1,z_2$, we can apply explicit Givens rotations to obtain a preferential tri-diagonal form for $H_{\theta,z_1,z_2}$, as described in~\cite{lam97}. In fact, the matrix then decomposes as the direct sum of two tridiagonal matrices, each of size approximately $q/2 \times q/2$. Finding eigenvalues for these tridiagonal matrices is an order $q^2$ operation, which is a fast computation in this application.

\section{Gap labelling: inverse slope 1, -1}

To draw the butterfly wings, we do not draw the spectral lines, but instead draw curves around the spectral gaps. 

The spectral gaps are conveniently labelled using a Diophantine equation, where integer parameters $t,s$ are fixed, and we select the r-th gap in spectrum $\theta = p/q$ using the formula
\[ r = t*p - s*q. \]
This indexing is known as gap labelling, as described in \cites{bel90,gold09,kam03,ypma07}, and the parameters $(t,s)$ are related to indices in K-theory, or equivalently, Chern numbers. Here, the point is  that at each spectral level $\theta = p/q$, the r-th gap does actually appear (if we also include the zero-length gap at spectral value 0), and we have a convenient way of marching through all the gaps.

Of course one must choose reasonable values of $t,s$ in order to give a reasonable $r$ value in the range of $1\ldots q-1$, since there are only $q-1$ gaps. Thus for integers $t >0$ we require \[0 \leq s \leq t-1 \]  and for $t <0$ we require \[ t \leq s \leq -1.\]
Note we are using a slightly different convention from the literature for our gap labelling, introducing a minus sign into the Diophantine equation
\[ r = t*p - s*q. \]
This allows us to keep indices $(t,s)$ with the same sign, which is just a matter of convenience in our algorithm.

As an example, we set $t=1$ which forces $s=0$, and we plot the gaps in first image shown in Figure~3. There are two curves in this first image. First, the curve on the left of the butterfly wing, which consists of the right endpoints of spectral intervals to the left of the gap. Second, the curve on the right of the butterfly, consisting of the left endpoints of the spectral intervals to the right of the gap. 

For the negative slope $t=-1$, we are forced to take $s=-1$, and we see the plot for the downward sloping butterfly wings, also shown in Figure~3.

\begin{figure}[ht] 
   \centering
   \includegraphics[width=2.8in]{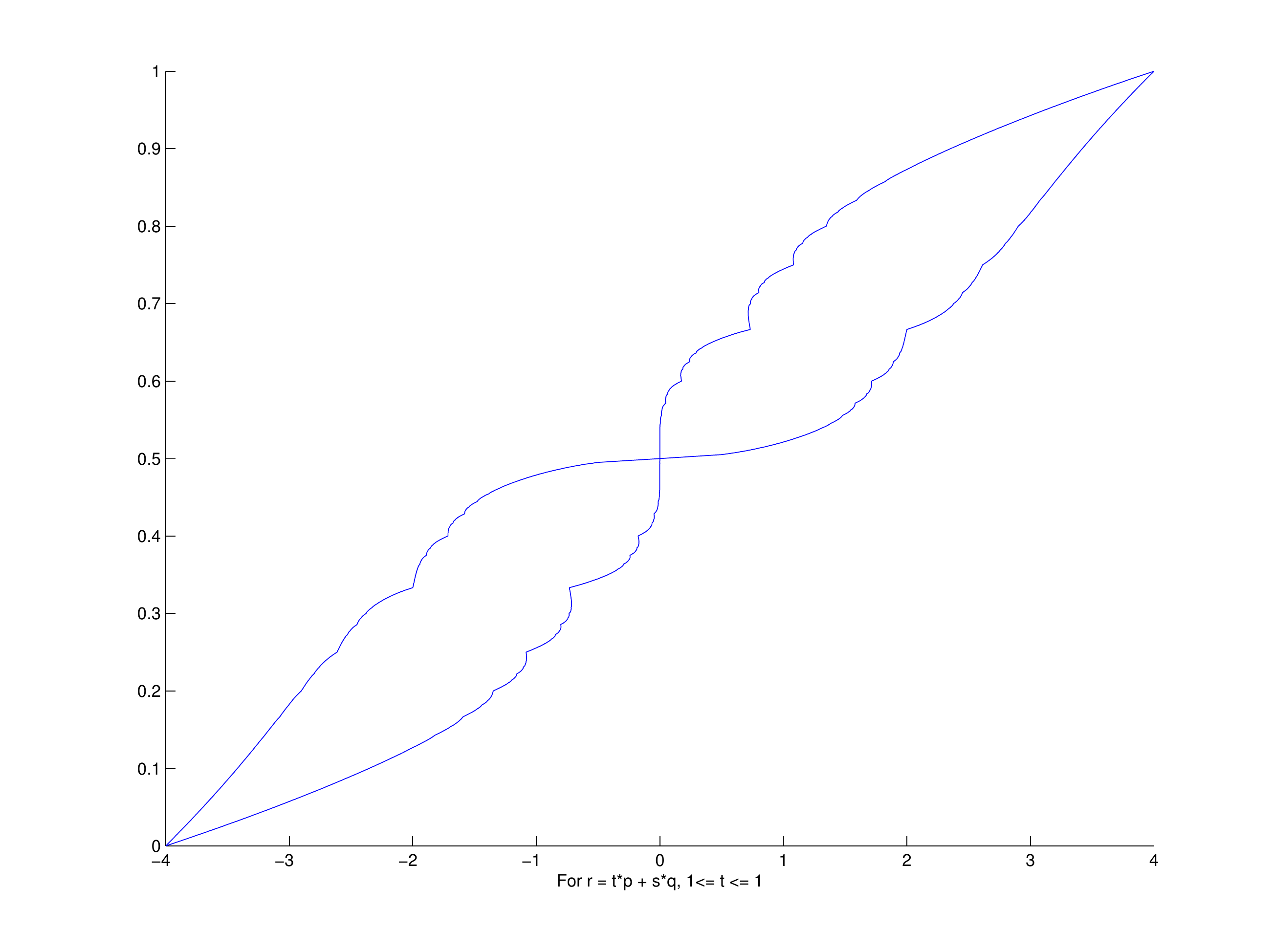} 
   \includegraphics[width=2.8in]{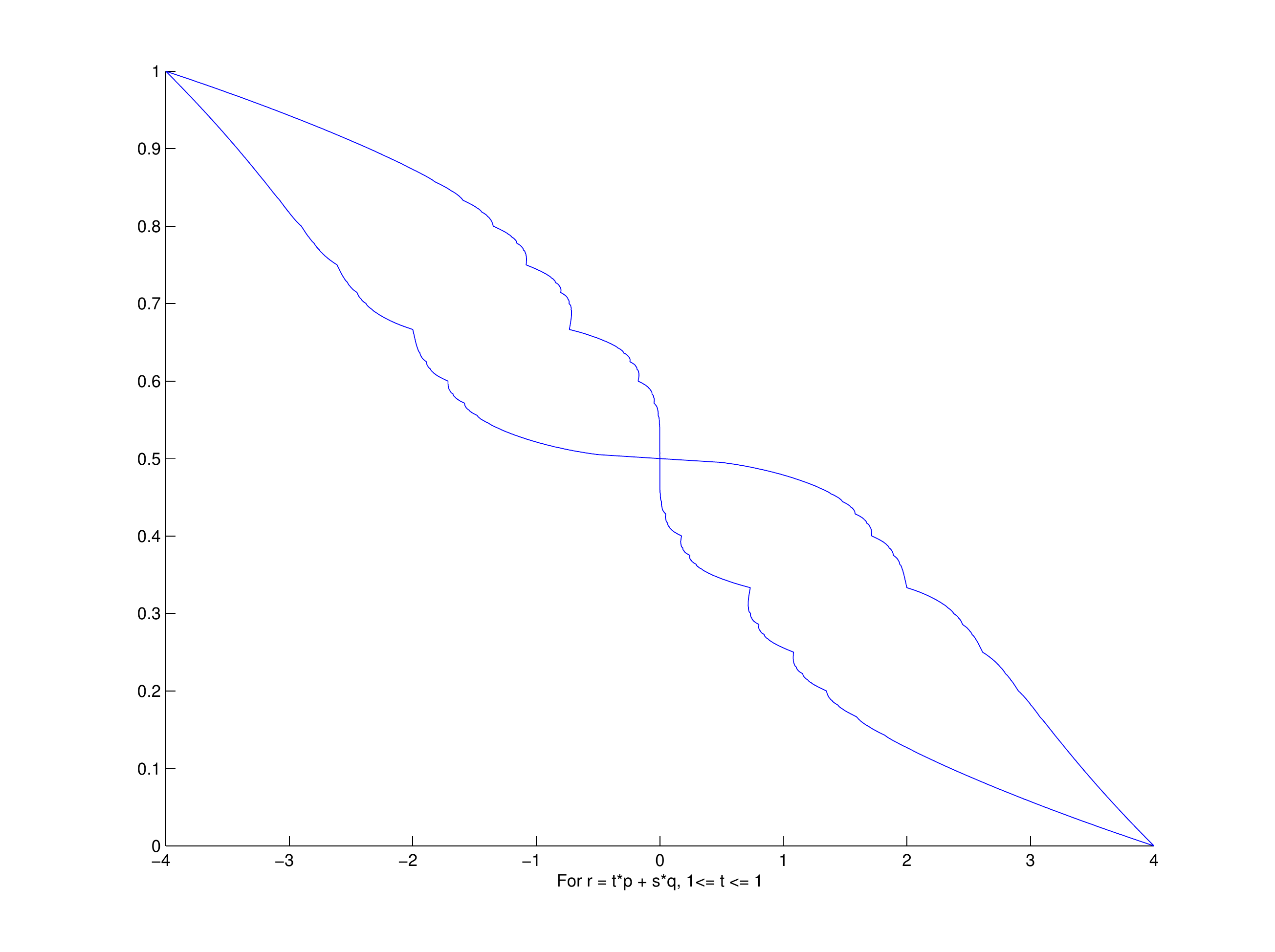} 
   \caption{Butterfly wings at slope 1 and -1.}
   \label{fig:BFs1}
\end{figure}

It is useful to refer to the parameter $t$ as the {\em inverse slope}, for the Diophantine equation
$r = t*p-s*q$ becomes, after division by $q$, a linear equation
\[ \frac{r}{q} = t*\frac{p}{q} - s, \mbox{ or }  x = t*\theta -s, \]
where $t$ indicates the slope between variables $x$ and $\theta$.

Notice something special happens at the wingtips, where $\theta = 0,1$. For these values of $\theta$, there is no spectral gap, as the spectrum consists of the single interval $[-4,4]$. So it does not make sense to identify the $r$-th gap there. From Figure~1, though, it is clear what should happen: the butterfly wing should close to a tip. To get this correct, we must include the endpoints of the interval $[-4,4]$ as the closing point for the wing.

We will see in the next sections that we need additional closures to complete the smaller butterfly wingtips. 

\section{Gap labelling: inverse slope 2}

We continue with the gap labelling, looking at gaps with label $r = 2*p -s*q$; that is, wings with inverse slope 2. For $s=0,1$, we get the wings shown in Figure~4, the first diagram covering a $\theta$ range of $[0,1/2]$ and the second covering range $[1/2,1]$. Note that we must close the wingtips at $\theta = 0, 1/2, 1,$ just as we did in the last section, selecting the min/max spectral endpoints, rather than a gap.

However, we also see a new phenomena. Again in Figure~4, on the left diagram, there is an extraneous segment connecting two distinct wings, at the vertical value $\theta = 1/3$. Similarly, in the the right diagram, there is a segment at $\theta = 2/3$. This is in fact a discontinuity in the gap labelling. At $\theta=1/3$, there is a jump in the gap. If we try to draw that jump line, we would get a line cutting across the larger butterfly wing drawn in the last section. (The wing with inverse slope $t=-1$.)

So, we should remove these discontinuities. Again, we want to close the wingtips, so at $\theta=1/3$ we choose the same spectral endpoint to close the wingtip. Specifically, when coming up from below, we pick a right endpoint of a spectral interval. When coming down from above, we pick a left spectral endpoint. This is exactly as we did with the endpoints $\theta = 0,1$.

With these discontinuities removed, we see proper wingtips, shown in Figure~5.

\begin{figure}[ht] 
   \centering
   \includegraphics[width=2.8in]{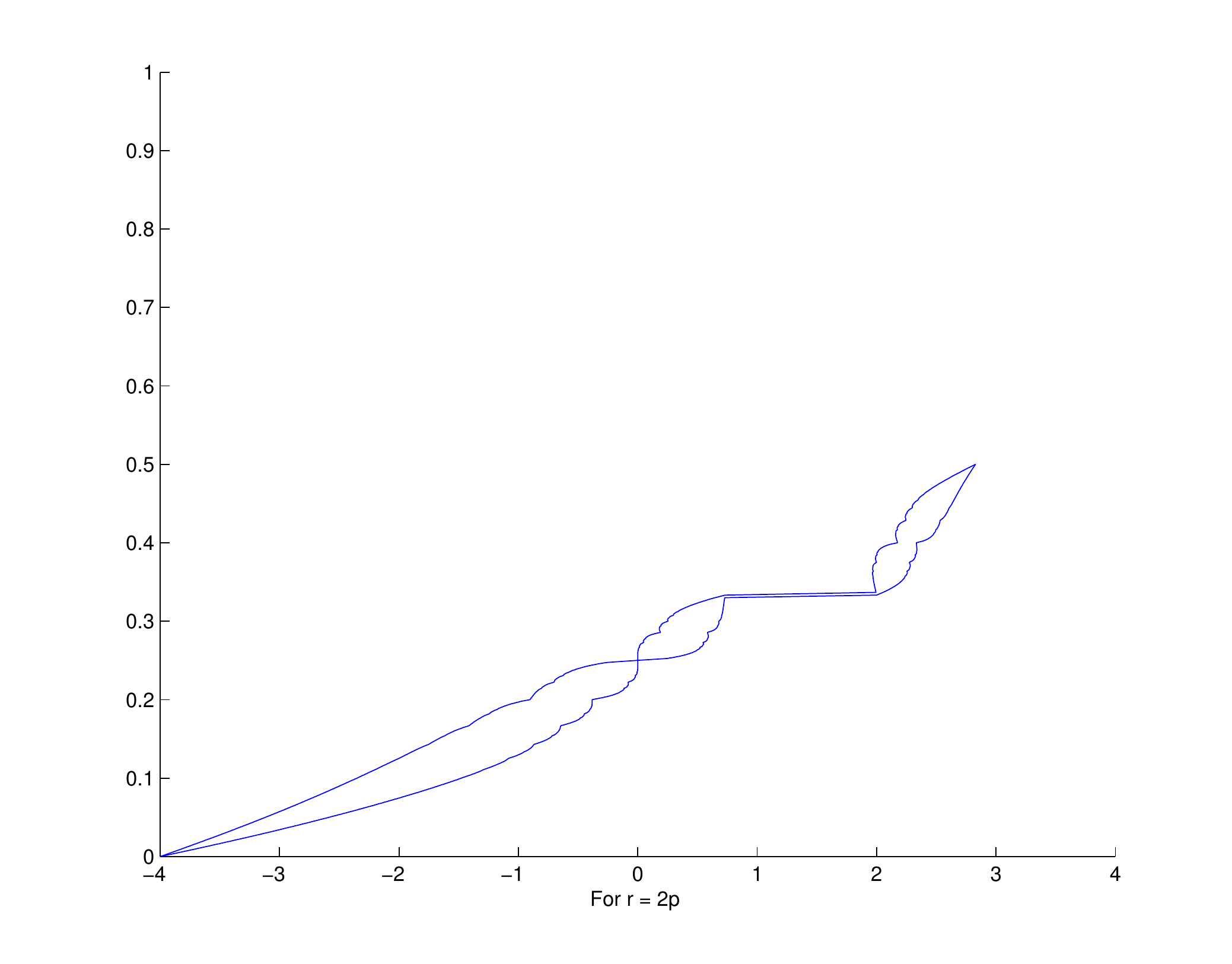} %
   \includegraphics[width=2.8in]{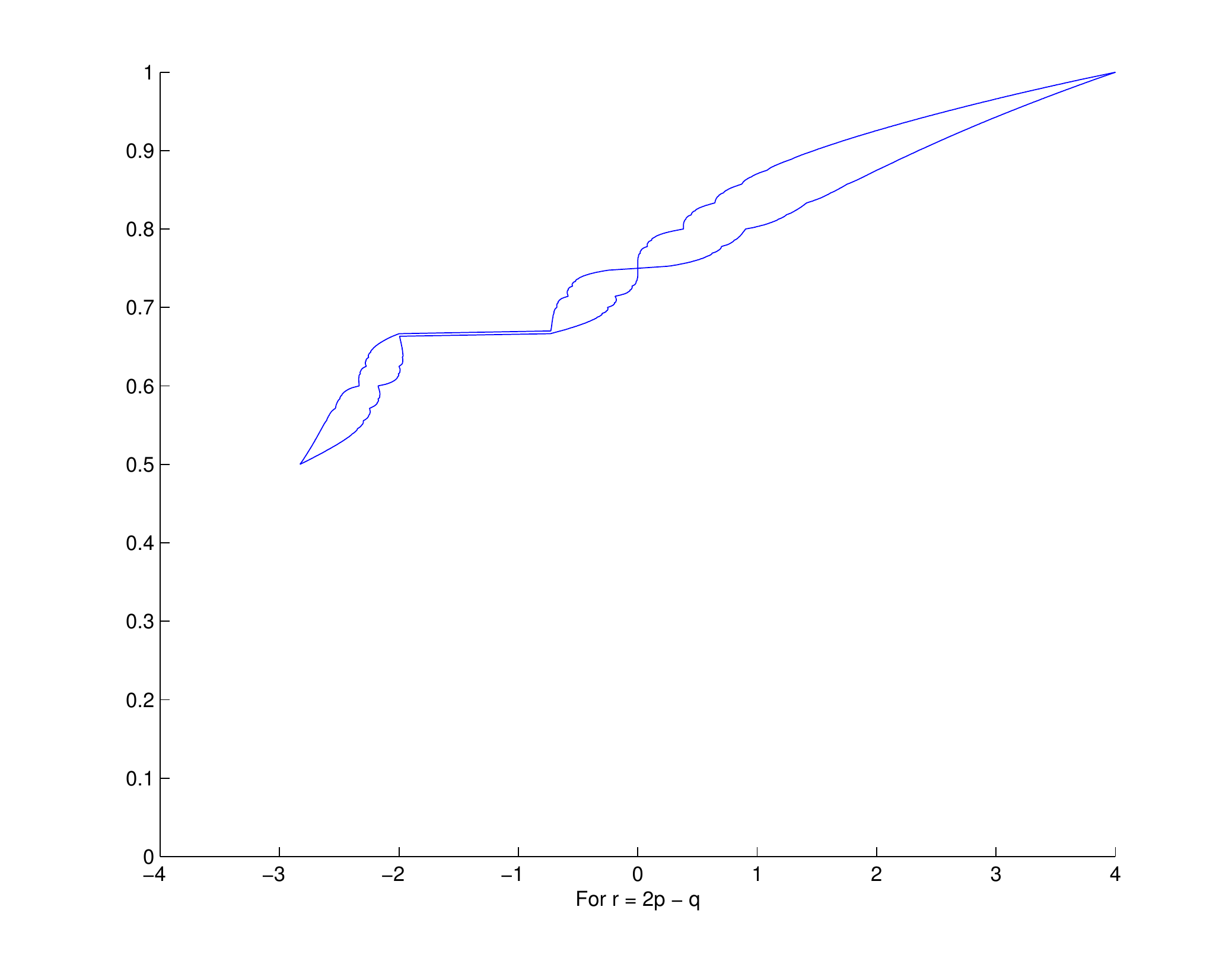} %
   \caption{Butterfly wings at inverse slope 2, $s=0,1$. Note the cutting segments at $\theta = 1/3, 2/3$.}
   \label{fig:BF2_jump}
\end{figure}

\begin{figure}[ht] 
   \centering
   \includegraphics[width=2.8in]{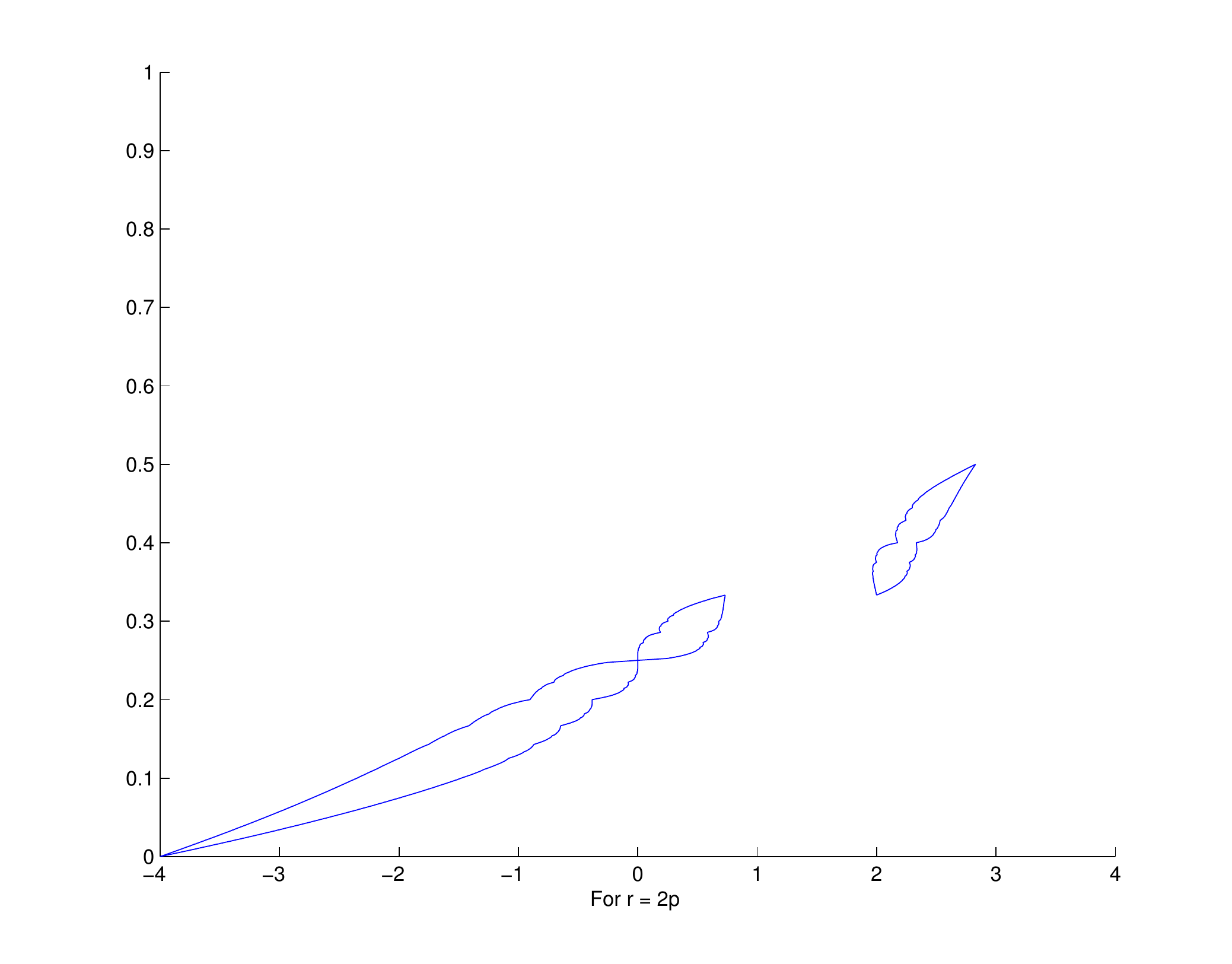} %
   \includegraphics[width=2.8in]{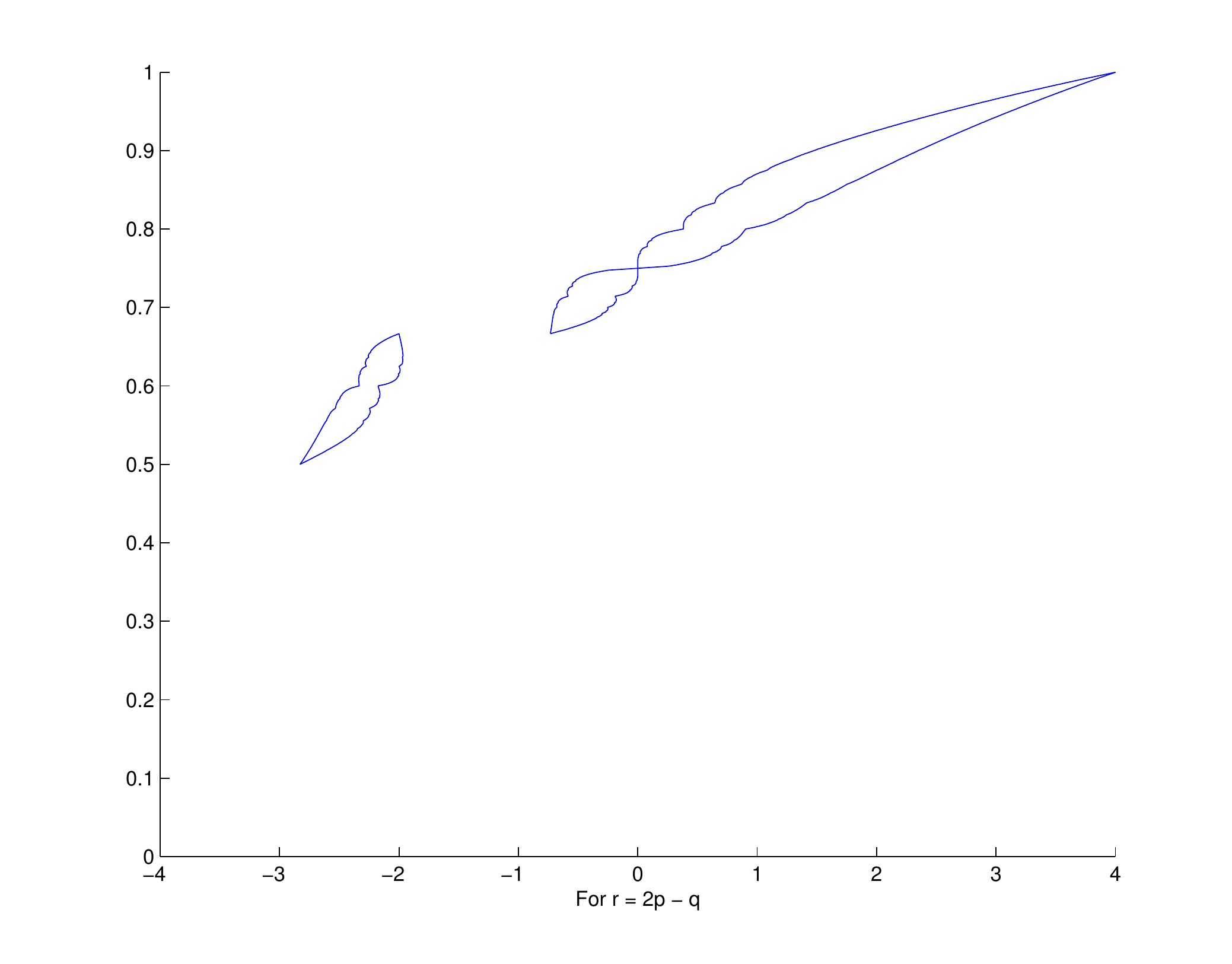} %
   \caption{Butterfly wings at inverse slope 2, $s=0,1$. Cutting segments removed.}
   \label{fig:BF2}
\end{figure}

\section{Gap labelling: inverse slope 3}

We continue the construction of butterfly wings, by labelling gaps with inverse slope $t=3$. In this case we have three choices for $s$, as $s=0,1,2$, giving the diagrams in Figures~6, 7, 8. Again, Figure~6 only covers the range of $\theta$ in $[0,1/3]$, Figure~7 covers the range $[1/3, 2/3]$ and Figure~8 the range $[2/3,1]$. We must close wingtips at the endpoints of these ranges (i.e. at $0, 1/3, 2/3, 1$), but we also notice cutting segments again. 

\begin{figure}[ht] 
   \centering
   \includegraphics[width=2.8in]{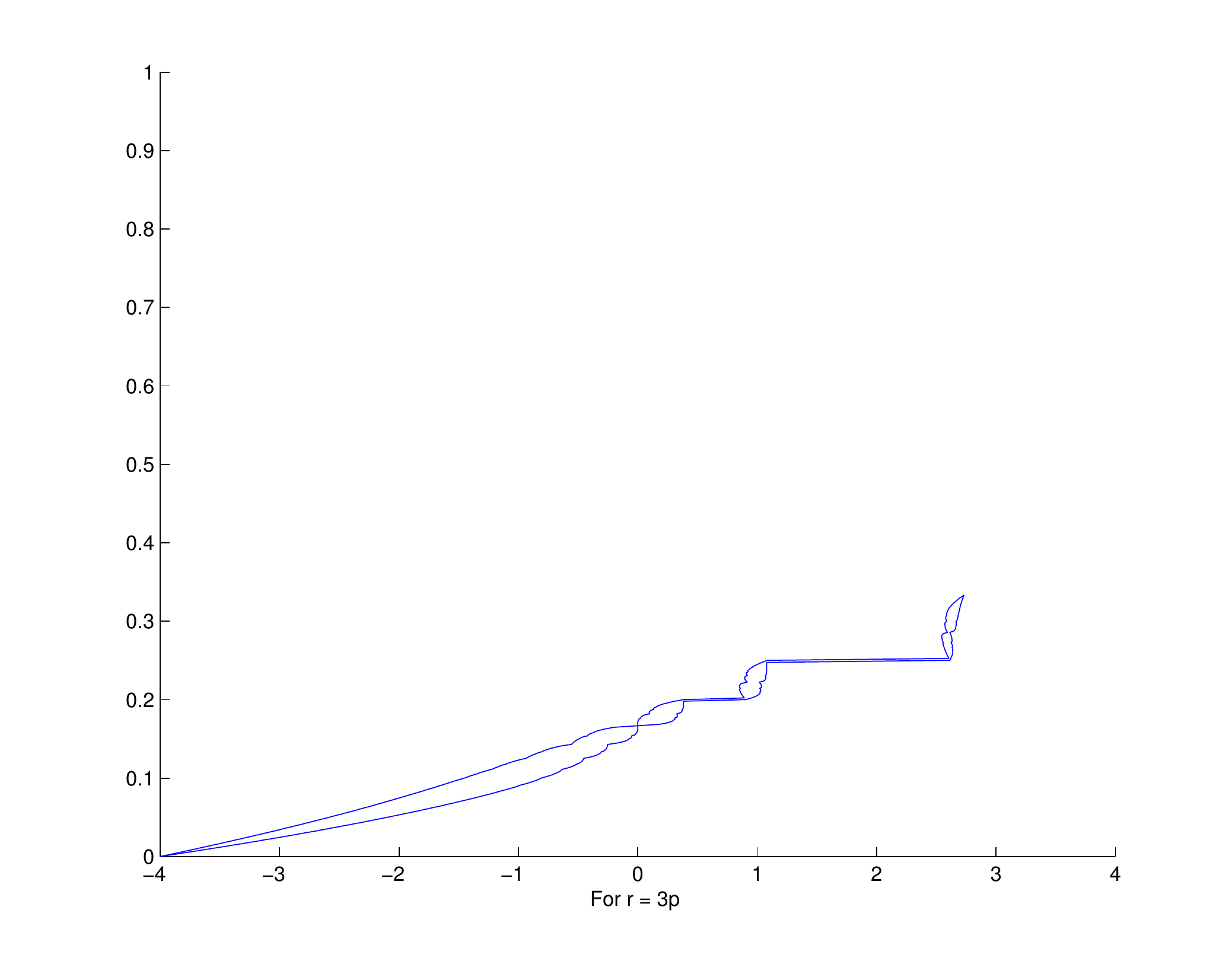} %
   \includegraphics[width=2.8in]{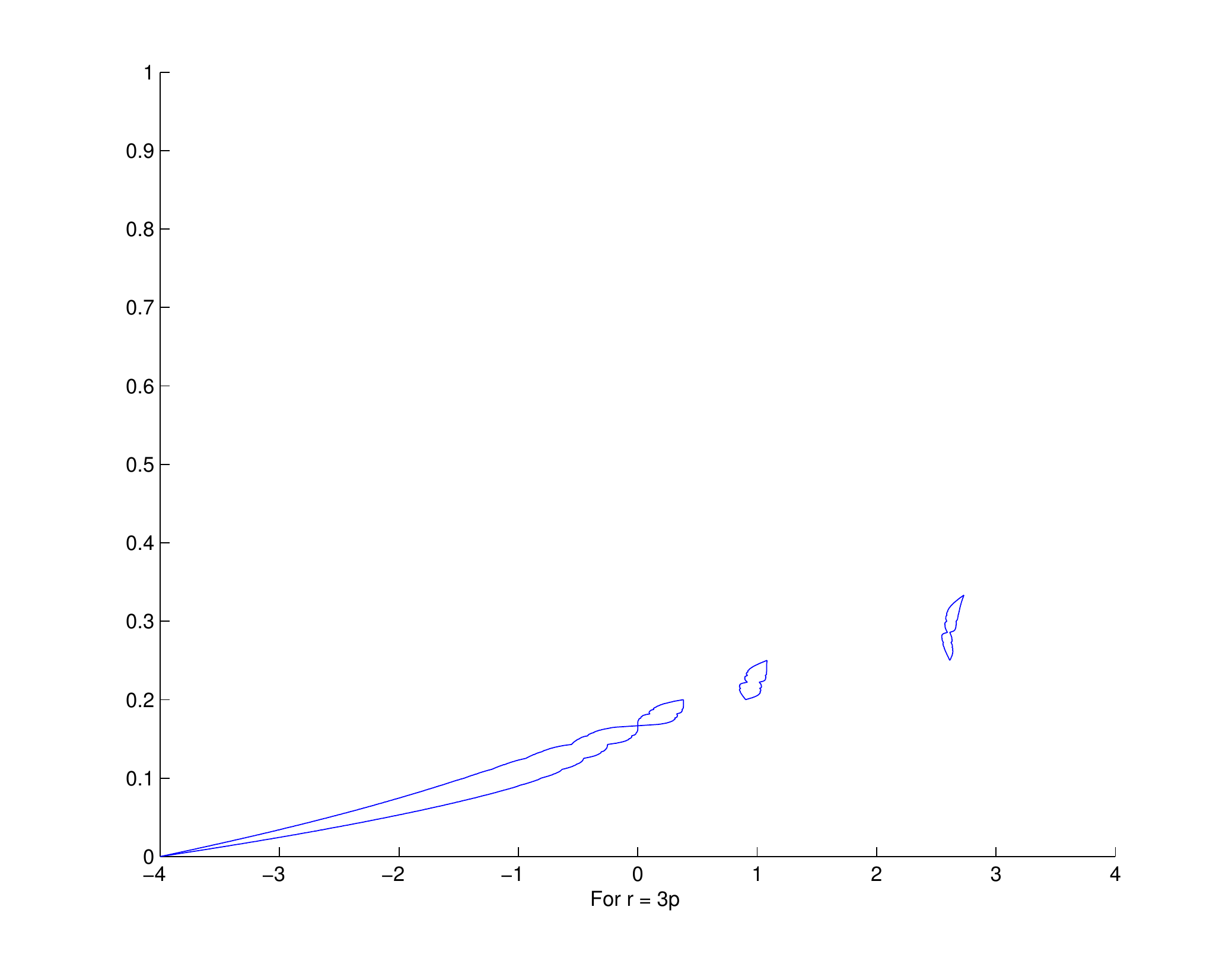} %
   \caption{Butterfly wings at inverse slope 3, $s=0$. Cutting segments on the left, at $\theta = 1/5, 1/4$, and removed on the right.}
   \label{fig:BF3a}
\end{figure}

\begin{figure}[ht] 
   \centering
   \includegraphics[width=2.8in]{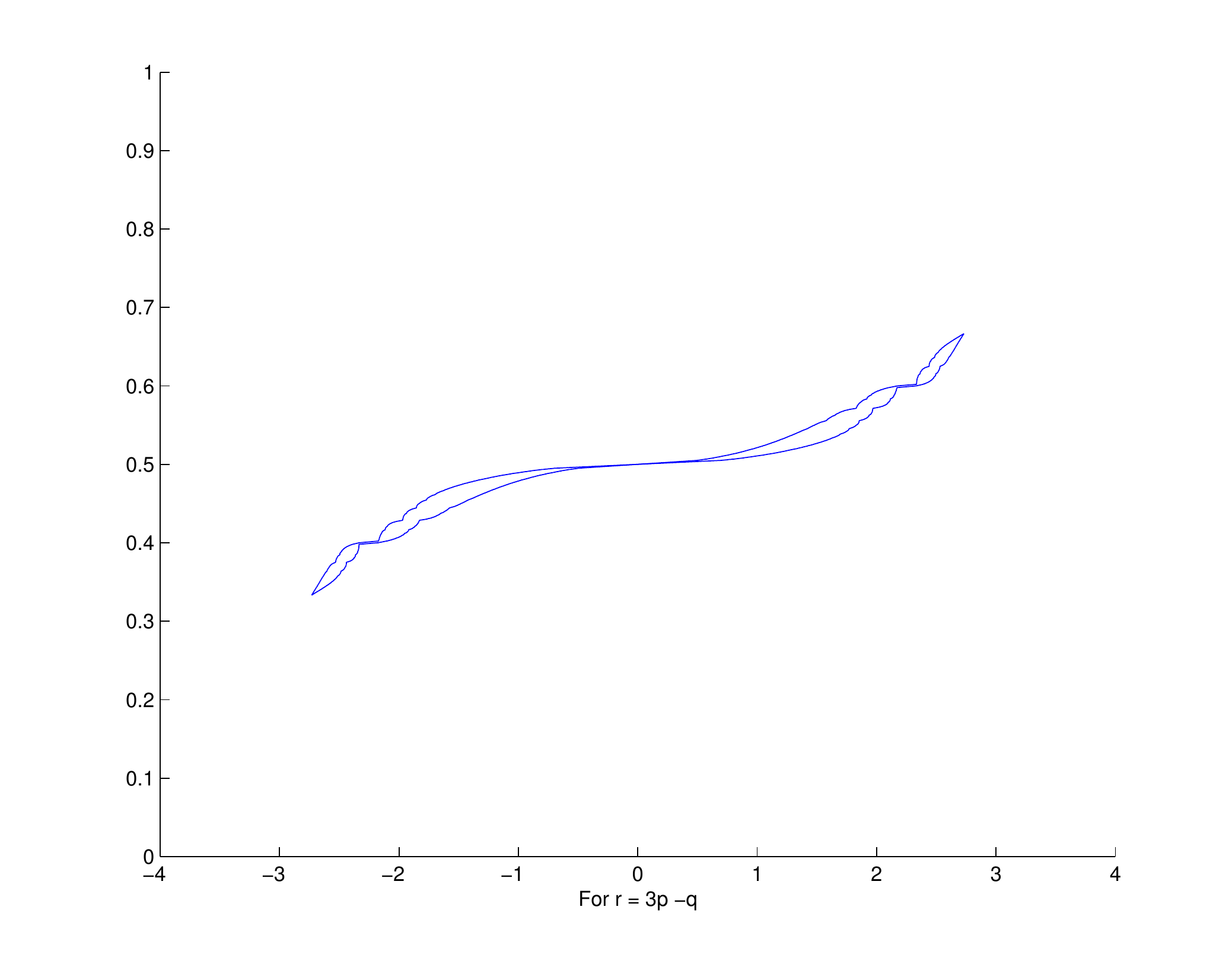} %
   \includegraphics[width=2.8in]{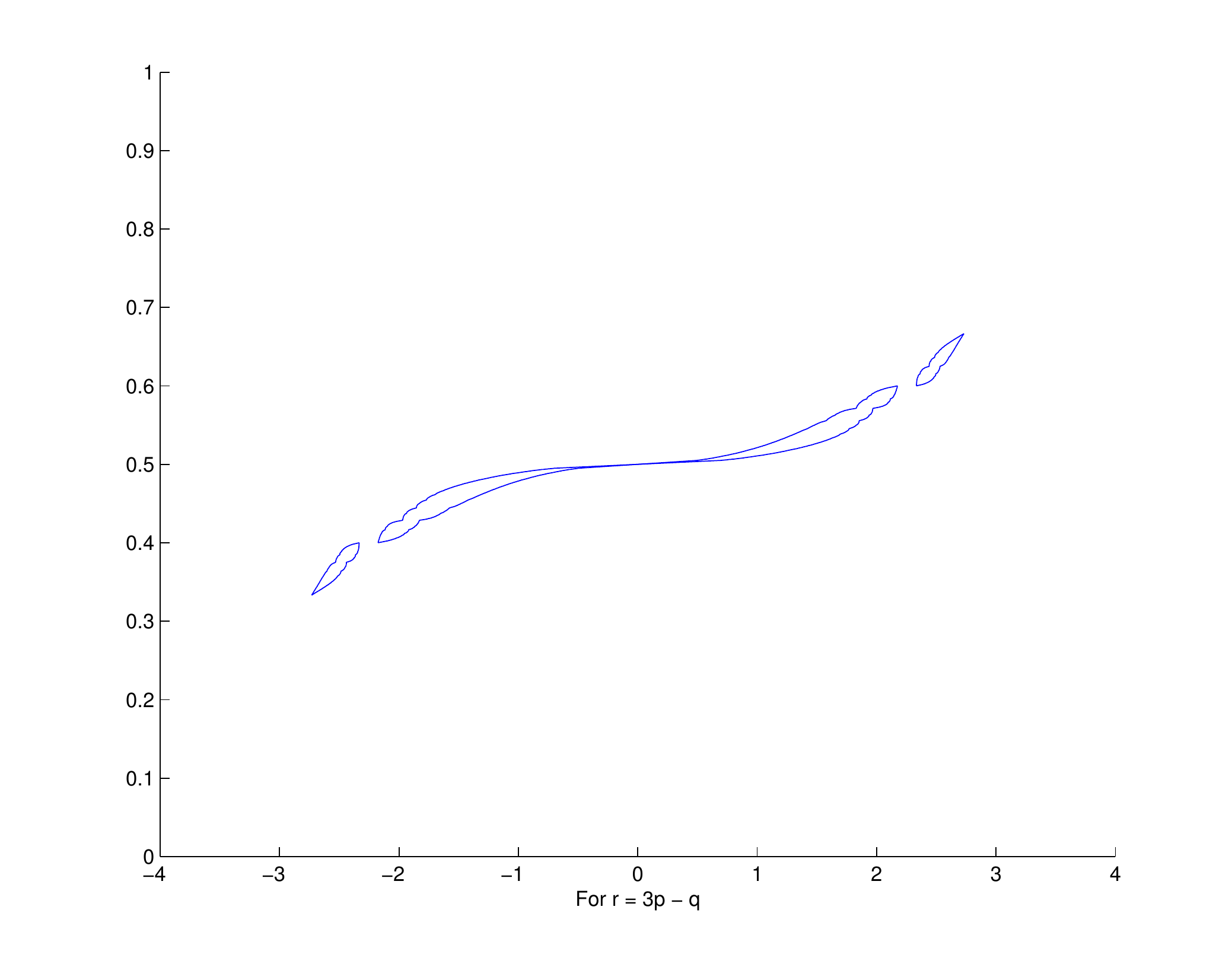} %
   \caption{Butterfly wings at inverse slope 3, $s=1$. Cutting segments on the left, at $\theta = 2/5, 3/5$, and removed on the right.}
   \label{fig:BF3b}
\end{figure}

\begin{figure}[ht] 
   \centering
   \includegraphics[width=2.8in]{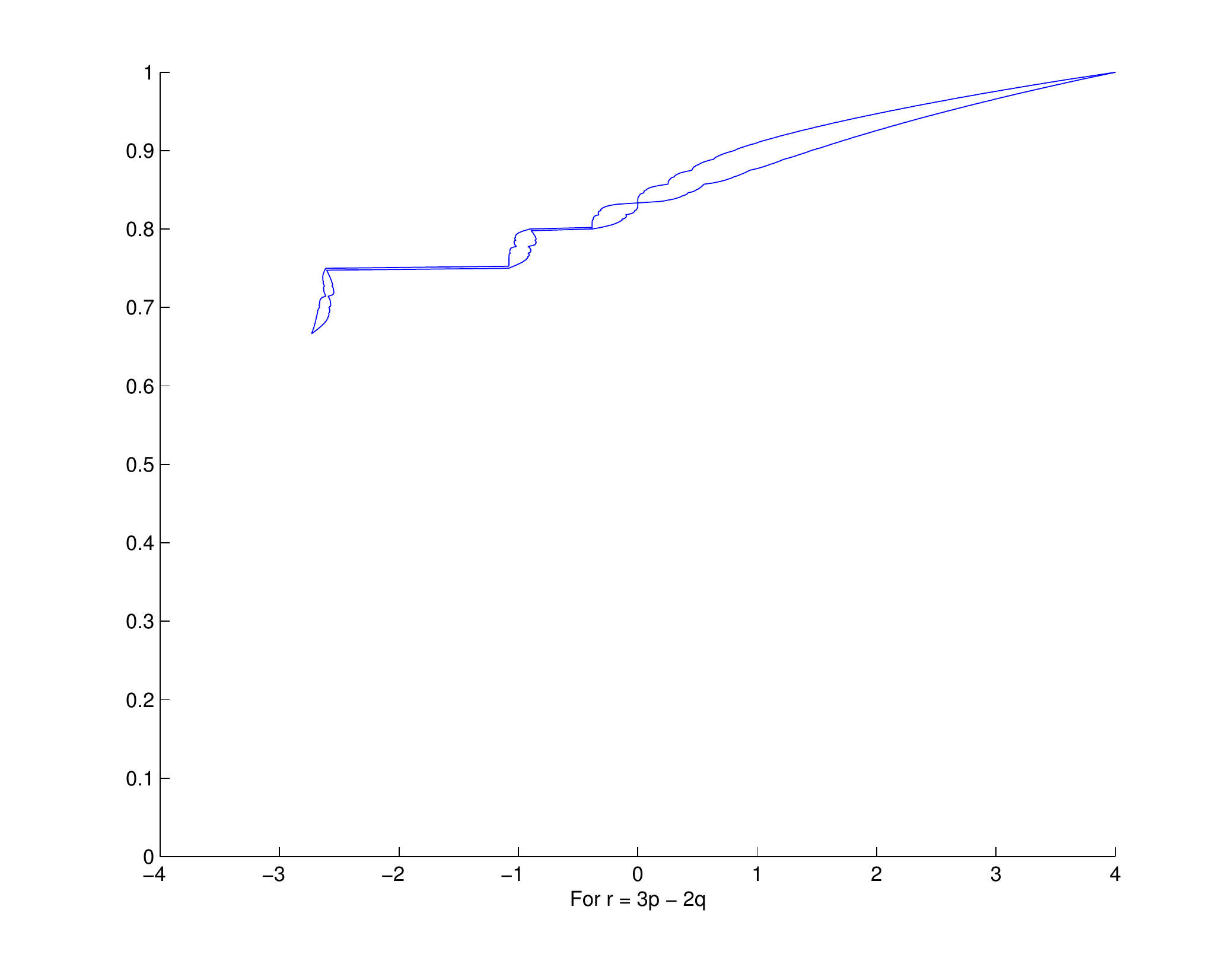} %
   \includegraphics[width=2.8in]{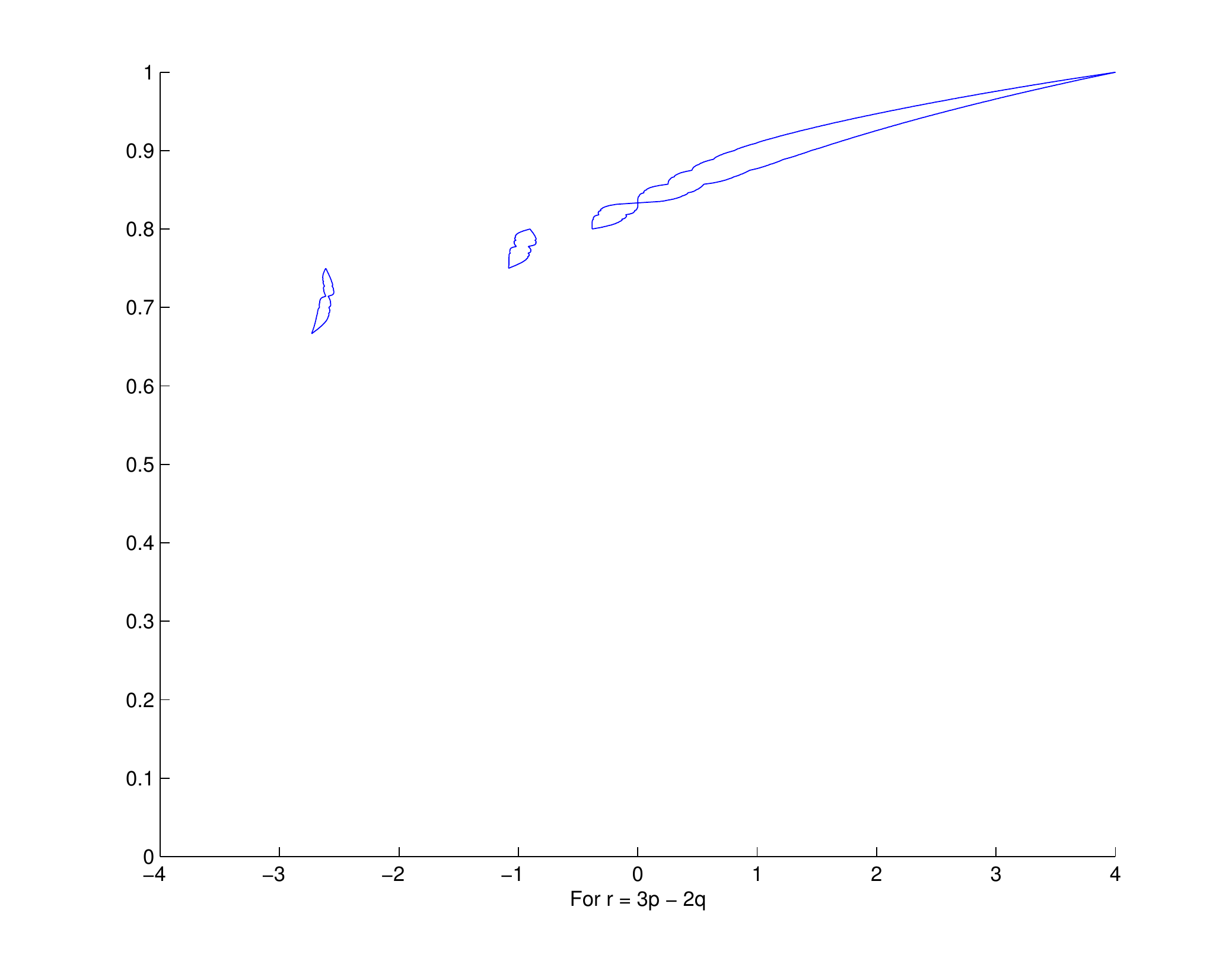} %
   \caption{Butterfly wings at inverse slope 3, $s=2$. Cutting segments on the left, at $\theta = 3/4, 4/5$, and removed on the right.}
   \label{fig:BF3c}
\end{figure}

In Figure~6, on the left, the cutting segments are clearly at $\theta = 1/5$ and $\theta=1/4$. Removing the cutting segments and closing the wingtips gives the right diagram in Figure~6.

In Figure~7, on the left, the cutting segments are not quite so clear, but a close examination shows there are segments  at $\theta = 2/5$ and $\theta=3/5$. Removing the cutting segments and closing the wingtips gives the right diagram in Figure~7.

In Figure~8, on the left, the cutting segments are  at $\theta = 3/4$ and $\theta=4/5$. Removing the cutting segments and closing the wingtips gives the right diagram in Figure~8.

\section{Gap labelling: inverse slope 4}

As we go to higher inverse slope $t$, it is harder to eyeball where the discontinuities should be.

For some ranges of $\theta$, the pattern is clear. For $t=4$, there will be four ranges of $\theta$, namely 
\[ [0, 1/4], \quad [1/4,1/2], \quad [1/2,3/4], \quad [3/4, 1]. \]
In the first interval, we expect wingtips (or discontinuities) at
\[ \theta = 0, 1/7, 1/6, 1/5, 1/4. \]
In the last interval, we expect wingtips (discontinuities) at
\[ \theta = 3/4, 4/5, 5/6, 6/7, 1, \]
which we can guess just from the patterns that occurred at slopes $t=2,3$.

\begin{figure}[ht] 
   \centering
   \includegraphics[width=4in]{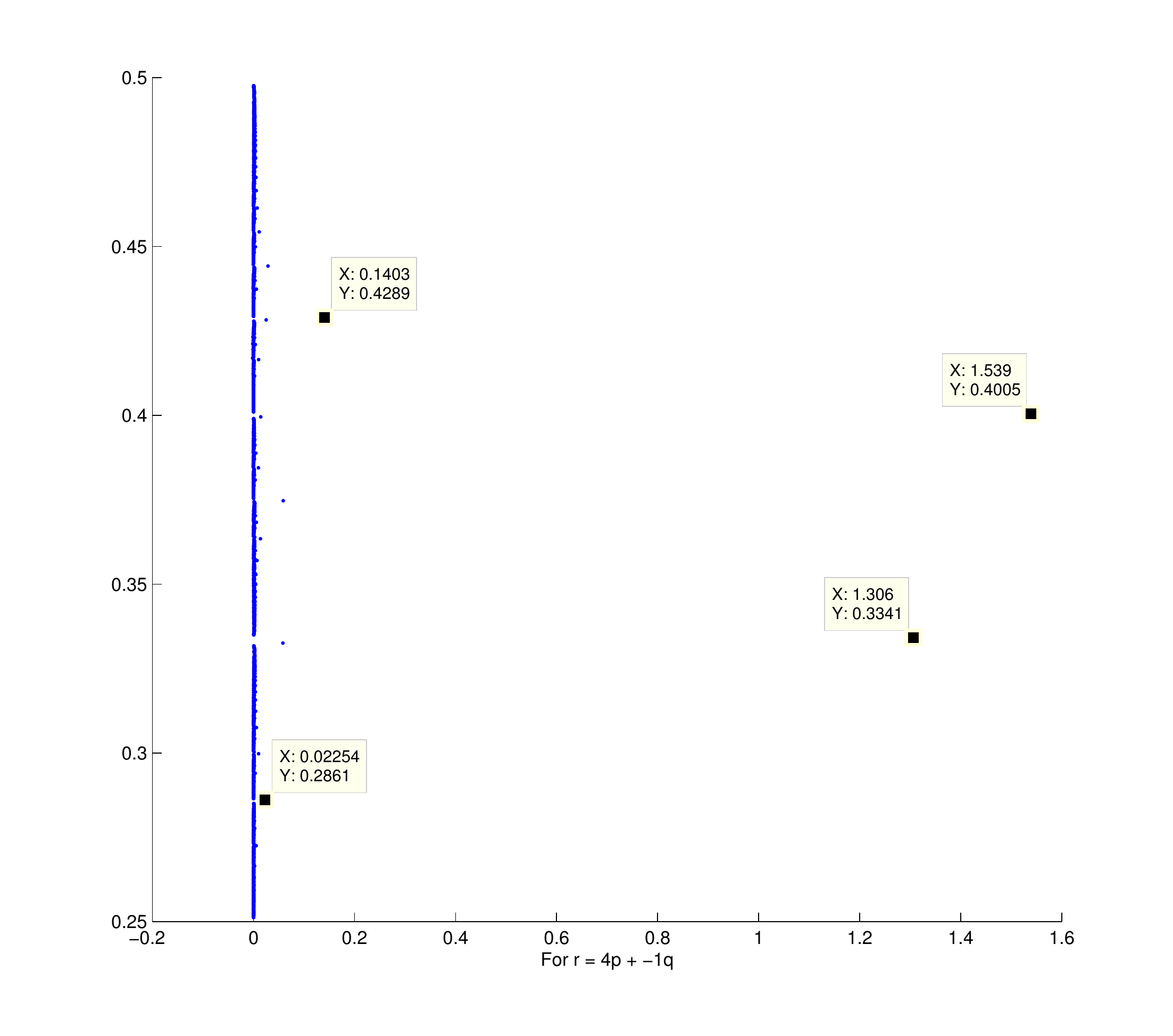} %
   \caption{A plot of jumps at inverse slope 4, big ones labeled.}
   \label{fig:gap_check}
\end{figure}

However, on the interval $[1/4,1/2]$, it is not so clear, so it helps to do some numerical experimentation. In Figure~9, we plot the jumps in the gap as $\theta$ increases. At the labelled points, we note significant jumps at the vertical  values $\theta=0.2861,0.3341,0.4005,0.4289$. 

Estimating in terms of fractional values, these jumps are at $\theta$ values  $2/7,1/3,2/5,3/7$. Significantly, there are four jump values -- whereas from the other examples, we might have only expected $3=t-1$ values. In Section~8, on the gap labelling conjecture, we show why there are four jump points. 

Using this jump data, we remove cutting lines and obtain the inverse slope 4 wings shown in Figure~10.

\begin{figure}[ht] 
   \centering
   \includegraphics[width=4in]{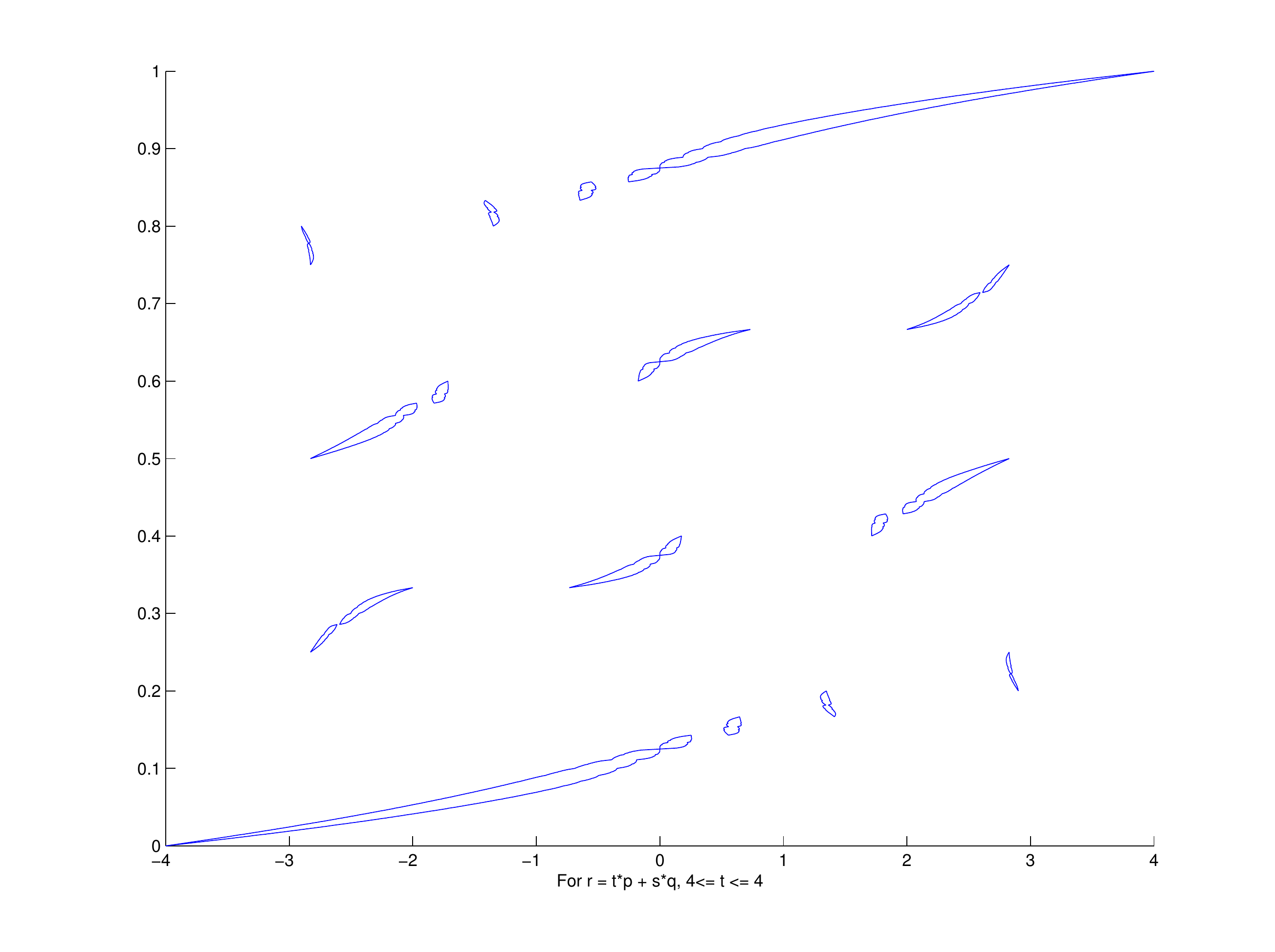} %
   \caption{Butterfly wings at inverse slope 4, $s=0,1,2,3$.}
   \label{fig:BFs4}
\end{figure}

\section{Sanity check: another jump at inverse slope 4?}

Looking carefully at Figure~9, it appears there is another significant gap jump at a $\theta$ value of about $.38$. Examining our numbers carefully in MATLAB, we find a value of $\theta = 0.3747$, which is about $\theta=3/8$.

Is there a jump here? No. This value corresponds to the spectral point zero, the centre of the butterfly pattern, where the inverse of the butterfly curves flattens out -- so we have large movement, even though there is no jump.

We note that $t*\theta - s$, at $\theta = 3/8$, gives the value $1/2$, the midpoint of the butterfly, or spectral value of zero.

\section{Gap labelling: inverse slope 5}

At inverse slope 5, it is very difficult to guestimate where the discontinuities should be. However, in the next section of this paper, we have a conjecture predicting the gap jumps for $\theta$ values of the form
\[ \theta = \frac{s+s_o}{t+t_o}, \mbox{ for any $0 < t < t_o$, $0<s\leq t$},\]
where we are following the gap curve $ r = t_o*p - s_o*q$, with $t_o > s_o \geq 0$.
For inverse slope 5, we set $t_o=5$ and $s_o = 0,1,2,3,4$.
We  use these fixed values, and the range of $t,s$ allowed above to pick out discontinuities, and then plot. Again, we must close the wingtips following the usual algorithm of matching endpoints of the appropriate spectral intervals. 

Figure~12 shows the collection of all the butterfly wings with inverse slope 5, all cutting lines removed, and wingtips closed, using the conjectured gap discontinuities. The image looks correct, suggesting that the conjecture works for inverse slope 5.

\begin{figure}[ht] 
   \centering
   \includegraphics[width=4in]{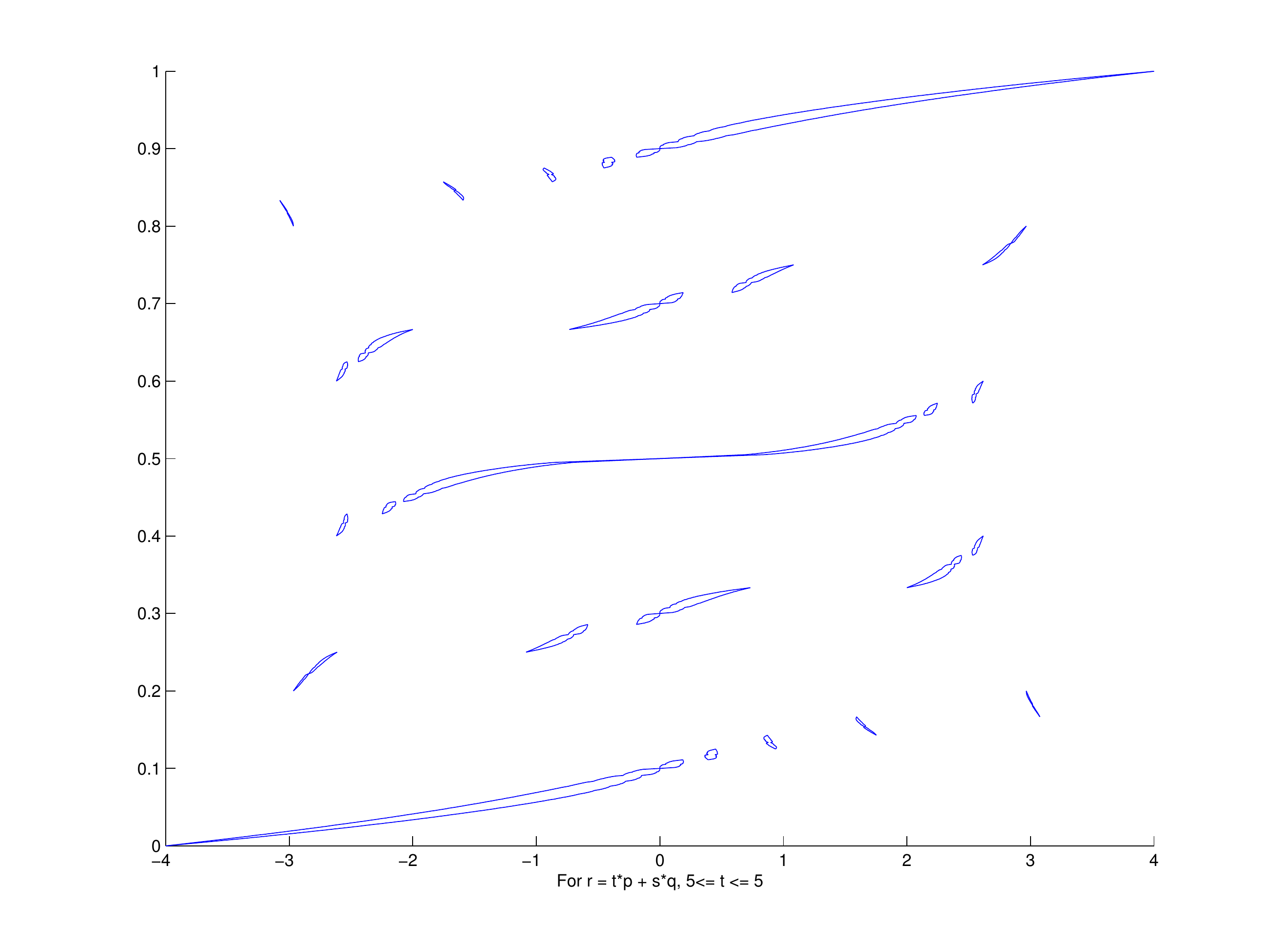} 
   \caption{Butterfly wings at inverse slope 5, $s=0,1,2,3,4$.}
   \label{fig:BFs5}
\end{figure}

\section{Gap labelling: discontinuity conjecture}

We observe the apparent discontinuities in the gap labelling can be determined by the following algorithm. 

{\em Algorithm}: For a given labelling $r = t_o*p - s_o q$, say with $t_o >0$, we plot the positive sloping line in $(x,\theta)$ given by the equation
\[ x = t_o * \theta - s_0, \]
and all the negative sloped lines in $(x,\theta)$ given by
\[ x = -t* \theta + s, \]
where $(t,s)$ are integers satisfying $0<s\leq t<t_o$. Note that $(-t,-s)$ is a valid index for a negative slope gap labelling. 

{\em Conjecture}: The intersection of all these lines, in the region $0<x<1$ give the gap discontinuities.

Visually, we can see the intersections. For instance, in the case of a line slope $t_o=2$, we see in Figure~12 the relevant lines for $s_0=0$ (on the left) and $s_o=1$ on the right, giving the intersections at $\theta=1/3$ and $\theta = 2/3$ respectively.

\begin{figure}[ht] 
   \centering
   \includegraphics[width=2.8in]{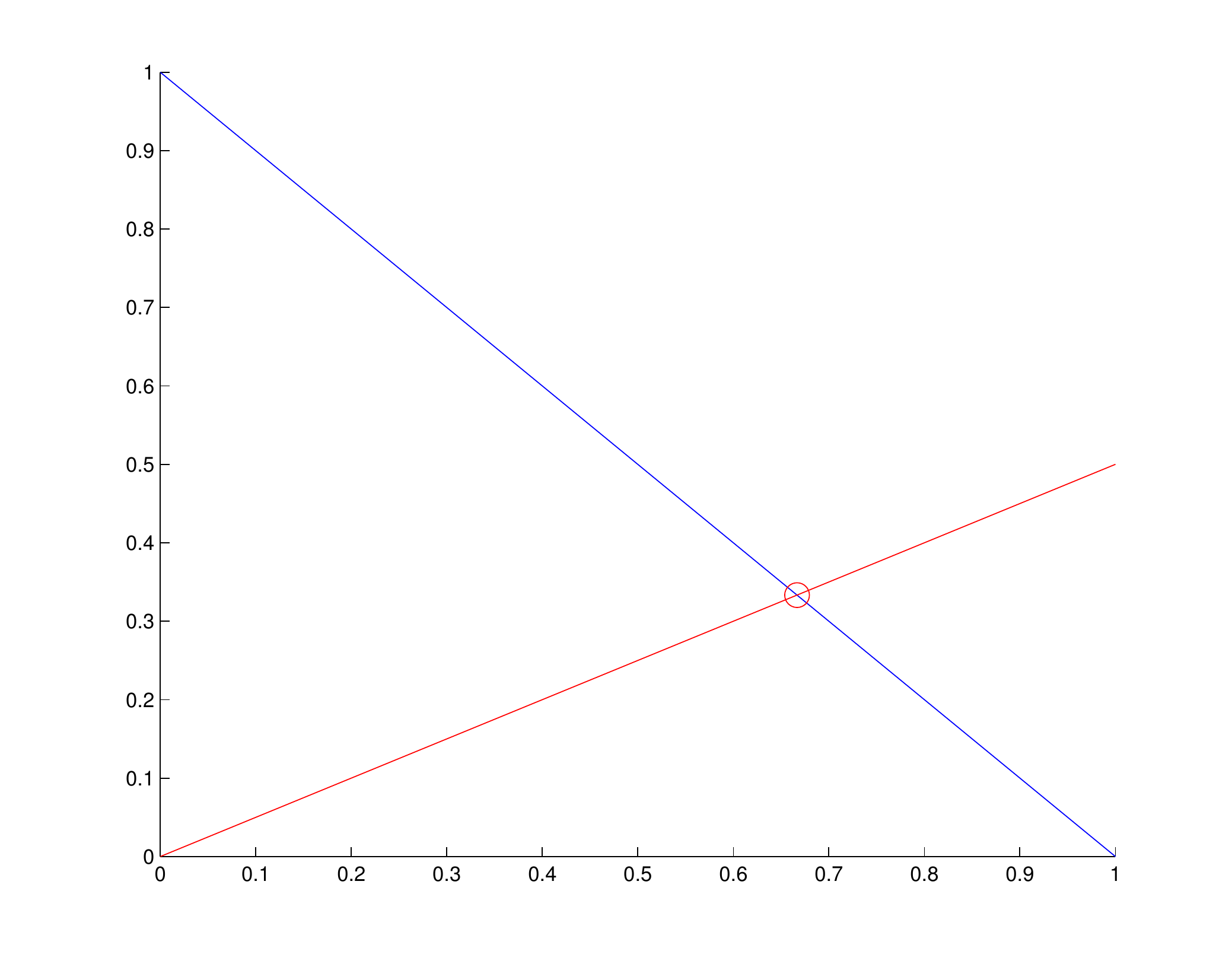} 
   \includegraphics[width=2.8in]{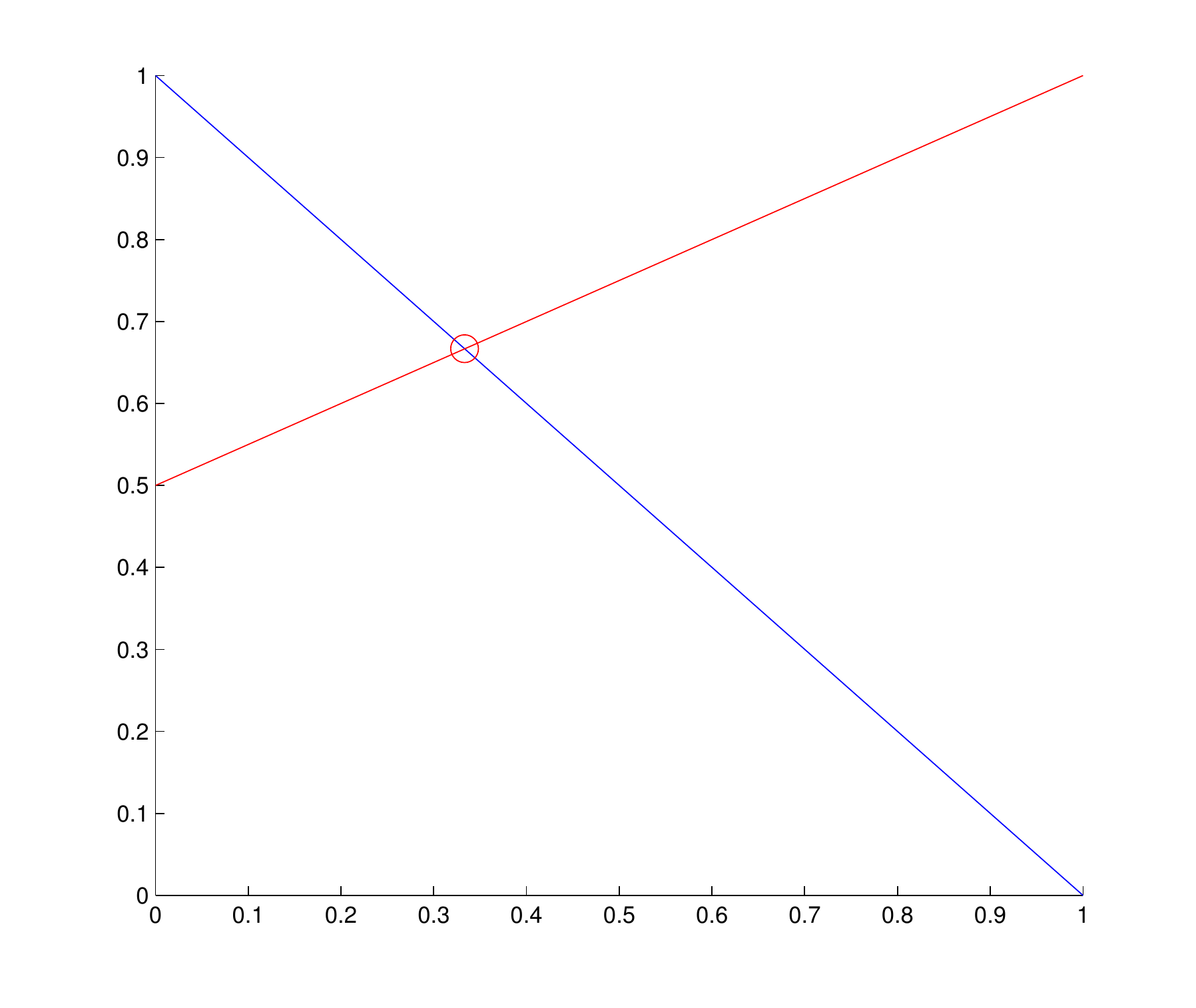} 
   \caption{Line intersections specify the discontinuity points for slope 2, $s=0,1$.
   On left, red line indexed as (2,0) has intersection at $\theta = 1/3$.
   On right, red line indexed as (2,1) has intersection at $\theta = 2/3$.}
   \label{fig:gaplines2}
\end{figure}

In the case of a line slope $t_o=3$, we see in Figure~13 the relevant lines for $s_0=0$, $s_o=1$, and $s_o=2$ has the intersections at $\theta=1/5,1/4$, $\theta=2/5,1/2,3/5$, $\theta=3/4,4/5$, respectively. The intersection at $\theta=1/2$ is not really a discontinuity, as the spectral intervals close here (spectral point zero). However, it does not hurt the algorithm to include this point as a jump point. 

\begin{figure}[ht] 
   \centering
   \includegraphics[width=1.8in]{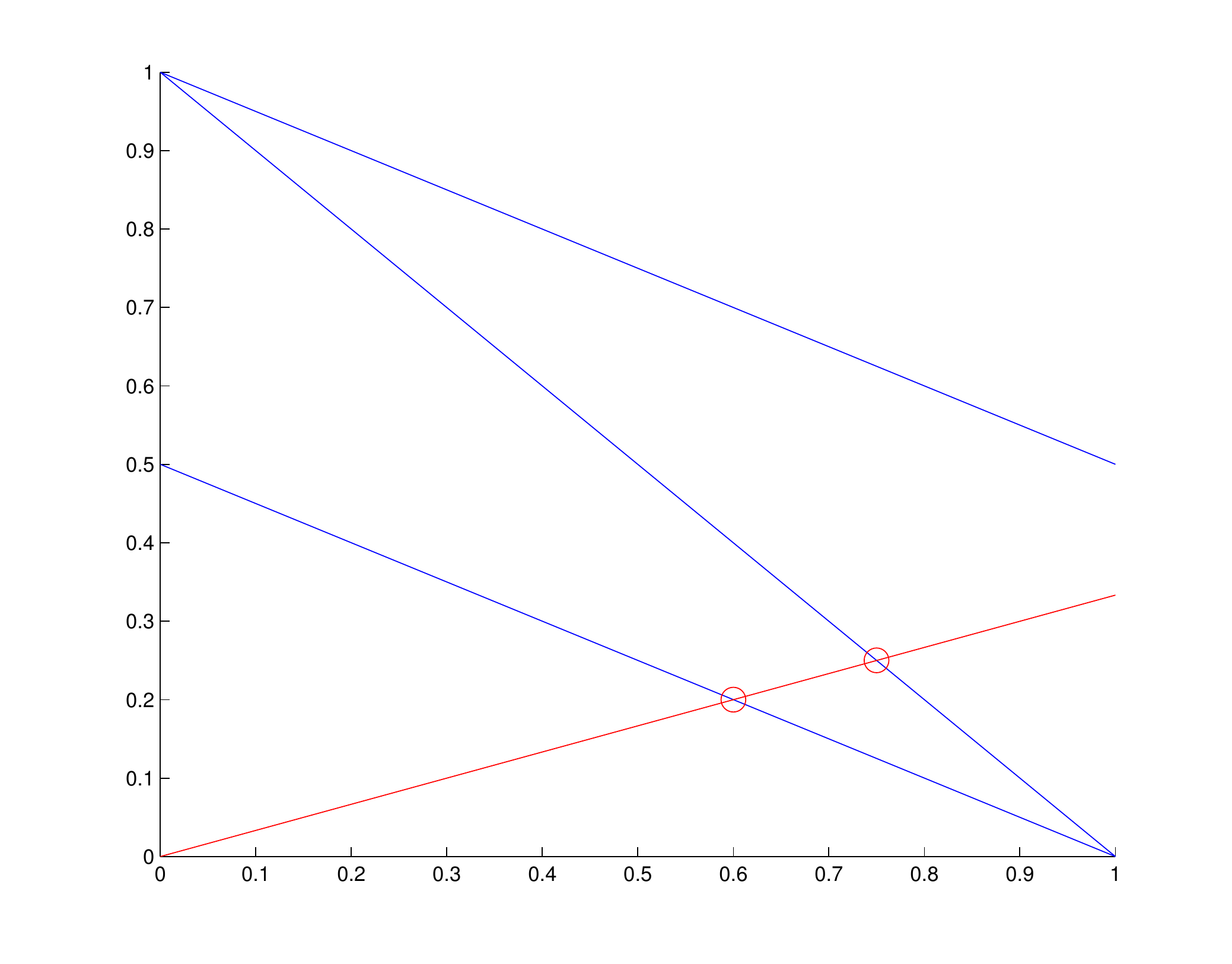} 
   \includegraphics[width=1.6in]{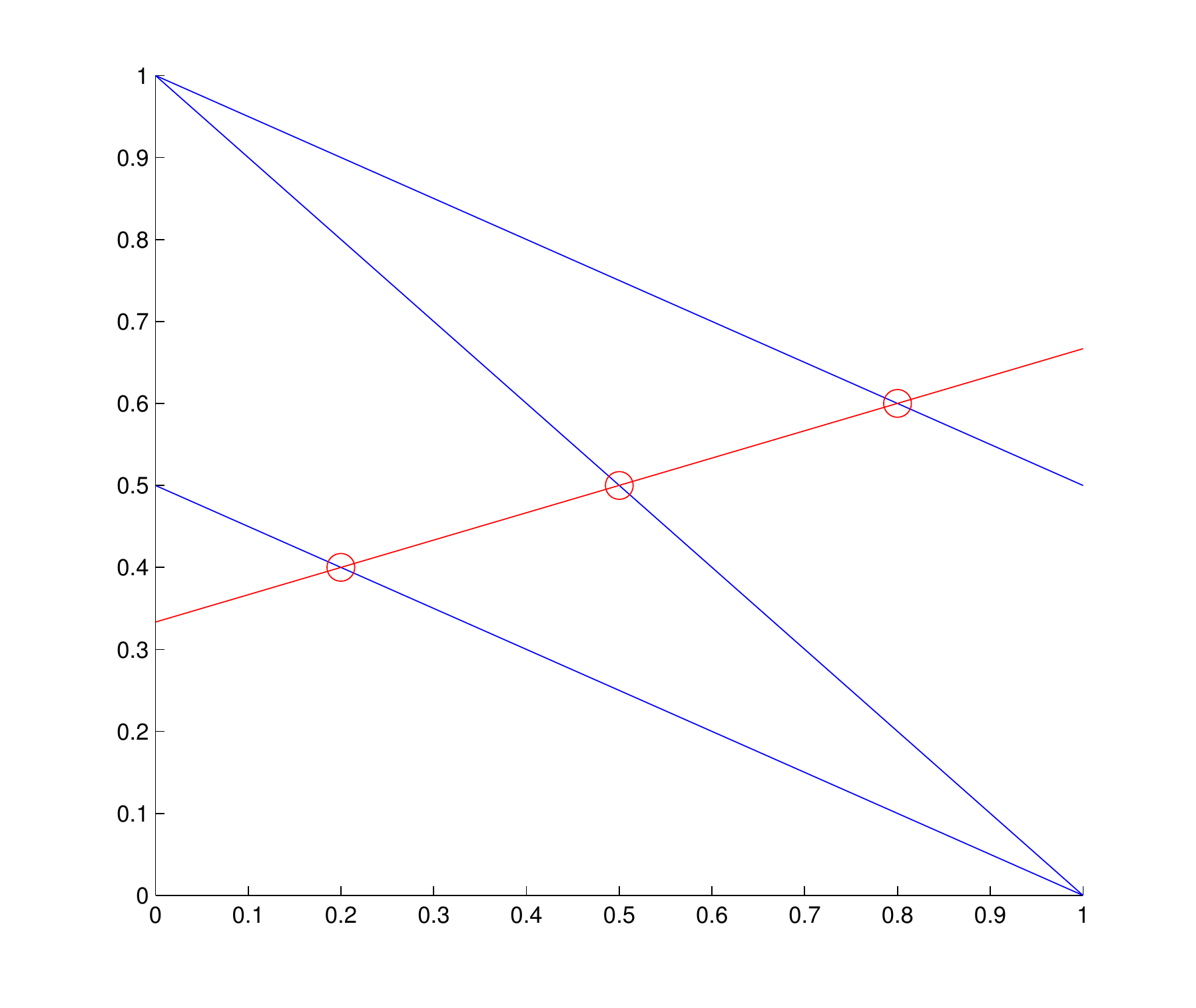} 
   \includegraphics[width=1.6in]{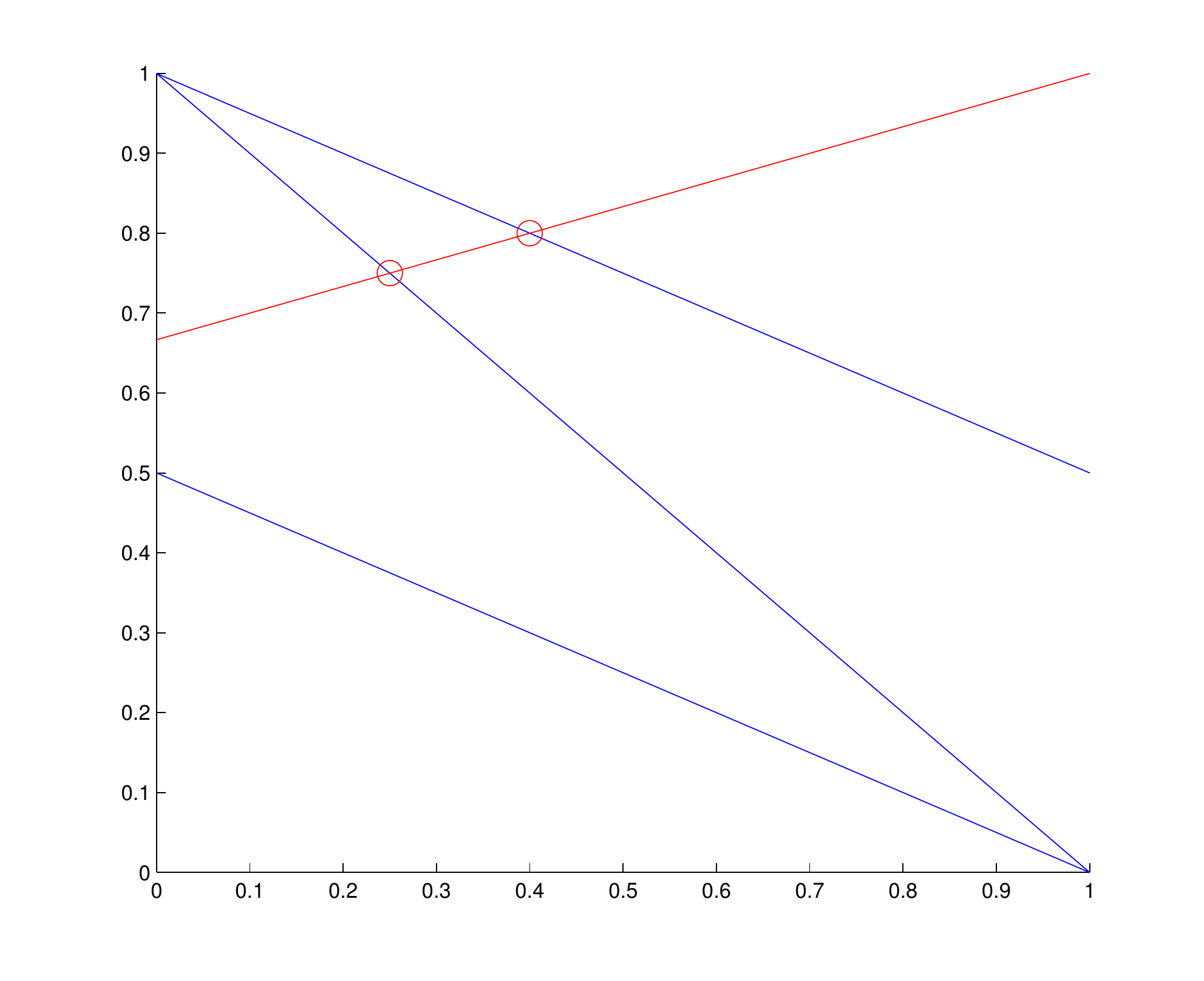} 
   \caption{Line intersections specify the discontinuity points for slope 3, $s=0,1,2$.
   Left, red line at (3,0) has intersections at $\theta = 1/5,1/4$.
   Centre, red line at (3,1) has intersections at $\theta = 2/5,1/2,3/5$.
   Right, red line at (3,2) has intersections at $\theta = 3/4,4/5$.}
   \label{fig:gaplines3}
\end{figure}

\begin{figure}[ht] 
   \centering
   \includegraphics[width=4in]{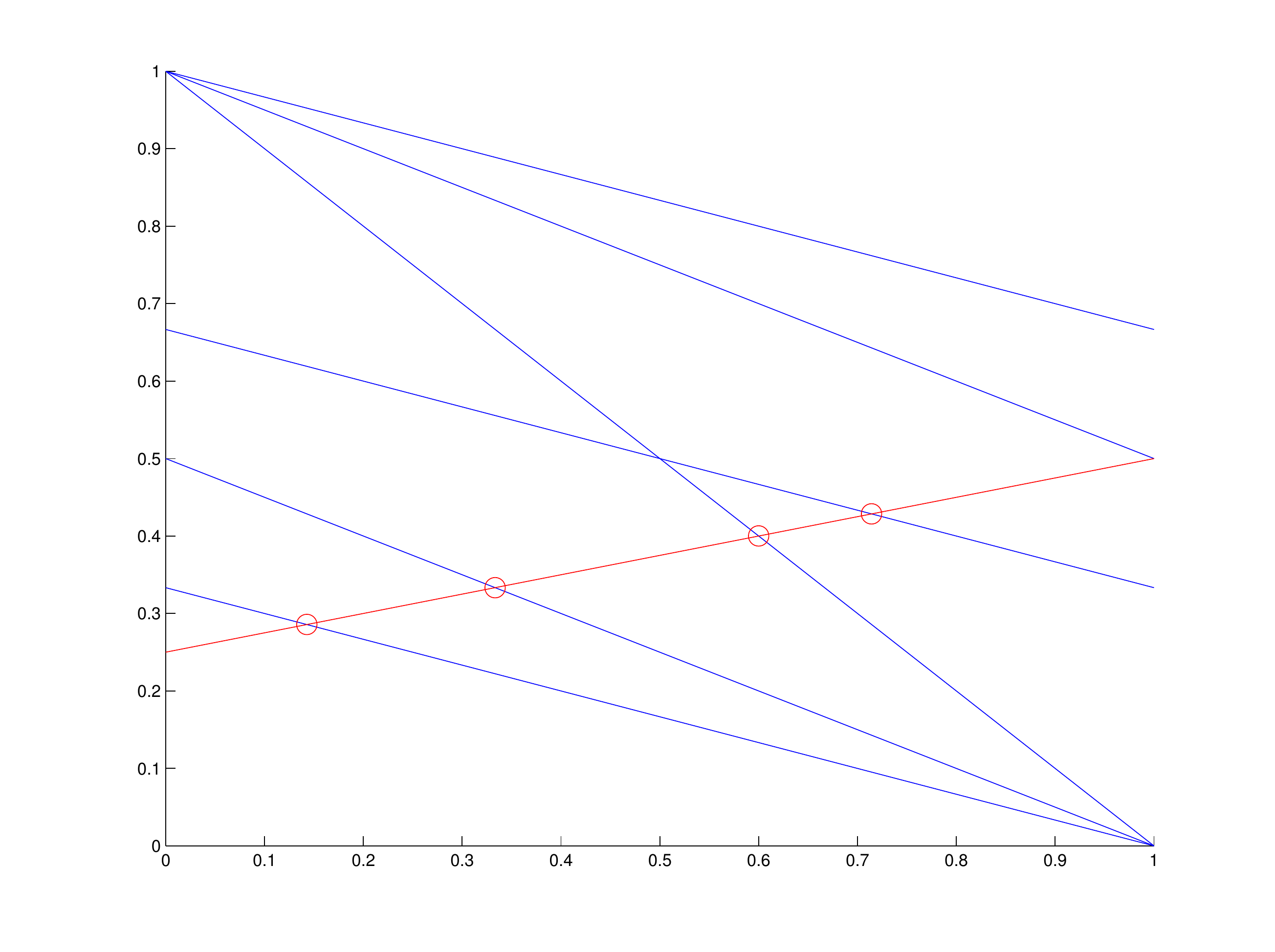} %
   \caption{Line intersections specify the discontinuity points.
   Red line indexed as (4,1) has four intersections, at $\theta = 2/7,1/3,2/5,3/7$.}
   \label{fig:gaplines4}
\end{figure}

The case of slope $t_o=4$ is particularly interesting, as this is the first level where the ``pattern'' of jump discontinuities is not immediately clear. Looking at Figure~14, for the line of index $(4,1)$, we see gap discontinuities predicted at values $\theta = 2/7,1/3,2/5,3/7$, just as we determined numerically in Section~6. This seems to confirm that the algorithm is promising (and perhaps correct).

It is useful to state specifically the conjecture about the gap discontinuities:

\begin{theorem}[Conjecture]
For gaps labelled by the Diophantine equation
\[ r = t_o * p - s_o *q, \mbox{ with $t_o > s_o \geq 0$}, \]
there are discontinuities in the gap labelling at values $\theta$ of the form
\[ \theta = \frac{s_o+s}{t_o+t}, \]
for all integers $t,s$ satisfying $0 <s\leq t<t_o$, and restricted to the open interval
\[ \theta = \frac{s_o+s}{t_o+t} \in \left( \frac{s_o}{t_o} , \frac{s_o+1}{t_o} \right). \]
Closing the gaps at the discontinuities gives the butterfly wingtips.
\end{theorem}

We do not have a proof of this result, only numerical evidence to confirm it. The predicted discontinuities also include the zero-length gaps in the spectral intervals, where two intervls touch at spectral value zero. Including these points in the plotting algorithm does not disturb our images.

As an example of implementing this theorem, with $(t_o,s_o) = (4,1)$, we expect discontinuities at
\[ \frac{1+1}{4+1}, \frac{1+1}{4+2}, \frac{1+1}{4+3}, \frac{1+2}{4+3} \]
which gives exactly the points $\theta = 2/5, 1/3, 2/7, 3/7$ as expected.

As the final numerical test of the conjecture, we note that the butterfly image in Figure~2 was computed using this algorithm for detecting discontinuities in the gap labelling, and we observe that Figure~2 is an excellent rendering of the butterflies apparent in the original spectral map in Figure~1. Figure~2 includes all gap labelling with all non-zero inverse slopes in the range $-10 \leq t \leq 10$. It looks correct, suggesting the algorithm to remove gap discontinuities is working correctly. 

In the appendix we include a figure with gap labelling for all non-zero inverse slopes in the range $-20 \leq t \leq 20$. Again, we see no obvious breakdown in the algorithm used to remove the gap discontinuities.

\section{Sanity check: zooming in on discontinuities}

As a check that our algorithm for predicting the discontinuities in the gap labelling is correct, we will examine closely the image at an unexpected discontinuity. 

As mentioned in Section~6, the discontinuity on line $(t_o,s_o)=(4,1)$ at $\theta = 2/7$ was unexpected. For this line, we were expecting only three discontinuities, as there are big jumps at the three values $\theta = 1/3, 2/5, 3/7$. The value $2/7$ doesn't quite fit the apparent pattern, and yet the conjecture says there should be a discontinuity there.

To test this, we plot side-by-side the butterfly image using spectral lines, and using the new algorithm for the drawing the labelled gaps, zoomed in at the possible discontinuity point $\theta = 2/7$. The two images are shown in Figure~15. Notice how the horizontal line exactly at the center (left image) has a gap -- so there really is a spectral gap at $\theta = 2/7$. The image on the right, showing the butterfly, accurately captures this gap -- that is, we see the main butterfly image here has a gap at the centre. So the algorithm successfully identified the gap discontinuity in this case.

It is significant, and worth noting, that the two wingtips at the centre of the diagram do not meet, and are slightly offset. It is an interesting phenomenon, accurately rendered both with spectral lines (left image in Figure~15) and with the butterflies (right image in Figure~15). 

\begin{figure}[ht] 
   \centering
   \includegraphics[width=2.8in]{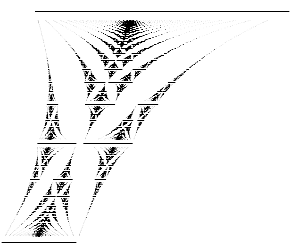} %
   \includegraphics[width=2.8in]{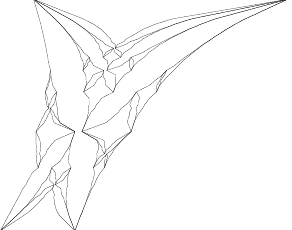} %
   \caption{Comparing the discontinuity at $\theta =2/7$, zoomed in.}
   \label{fig:zoom27}
\end{figure}

\section{Sanity check: other coupling constants}

It is a useful check on the algorithm to verify that the butterfly structures maintain the same basic form under small modifications of the setup of the problem. The Bloch electron model includes a physical parameter called the coupling constant $\lambda$ which is generically set equal to 2. For the general problem, one can consider self-adjoint operators of the form
\[ h_\theta = u + u^* + \frac{\lambda}{2}(v+v^*), \]
where $\lambda$ is a positive real parameter. It is known that the gap structure for the corresponding spectra remains the same~\cite{bel82}. 

We run our code as before, only changing the coupling parameter $\lambda$ and see that the basic form of the result looks the same. FIgure~7 shows three results, for three different values of the coupling parameter. Note the width of the result varies with $\lambda$, and in fact the width is proportional to $2 + \lambda$. Other than the width differences, the butterfly structures look very similar, and we do not observe any breakdown in the algorithm that predicts the gap discontinuities. Again this is evidence of a good conjecture.

\begin{figure}[ht] 
   \centering
   \includegraphics[width=1.63in]{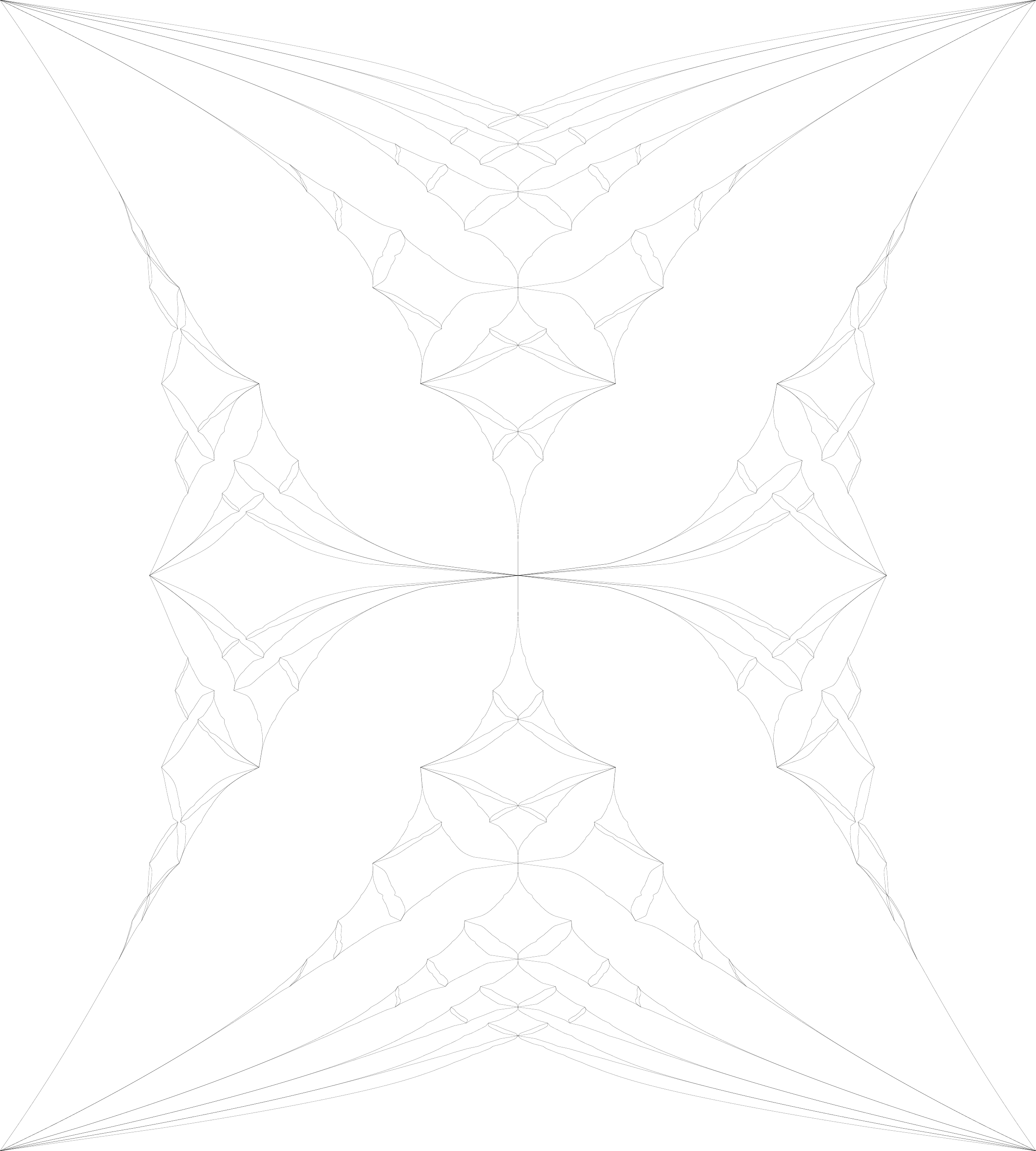} 
   \includegraphics[width=1.80in]{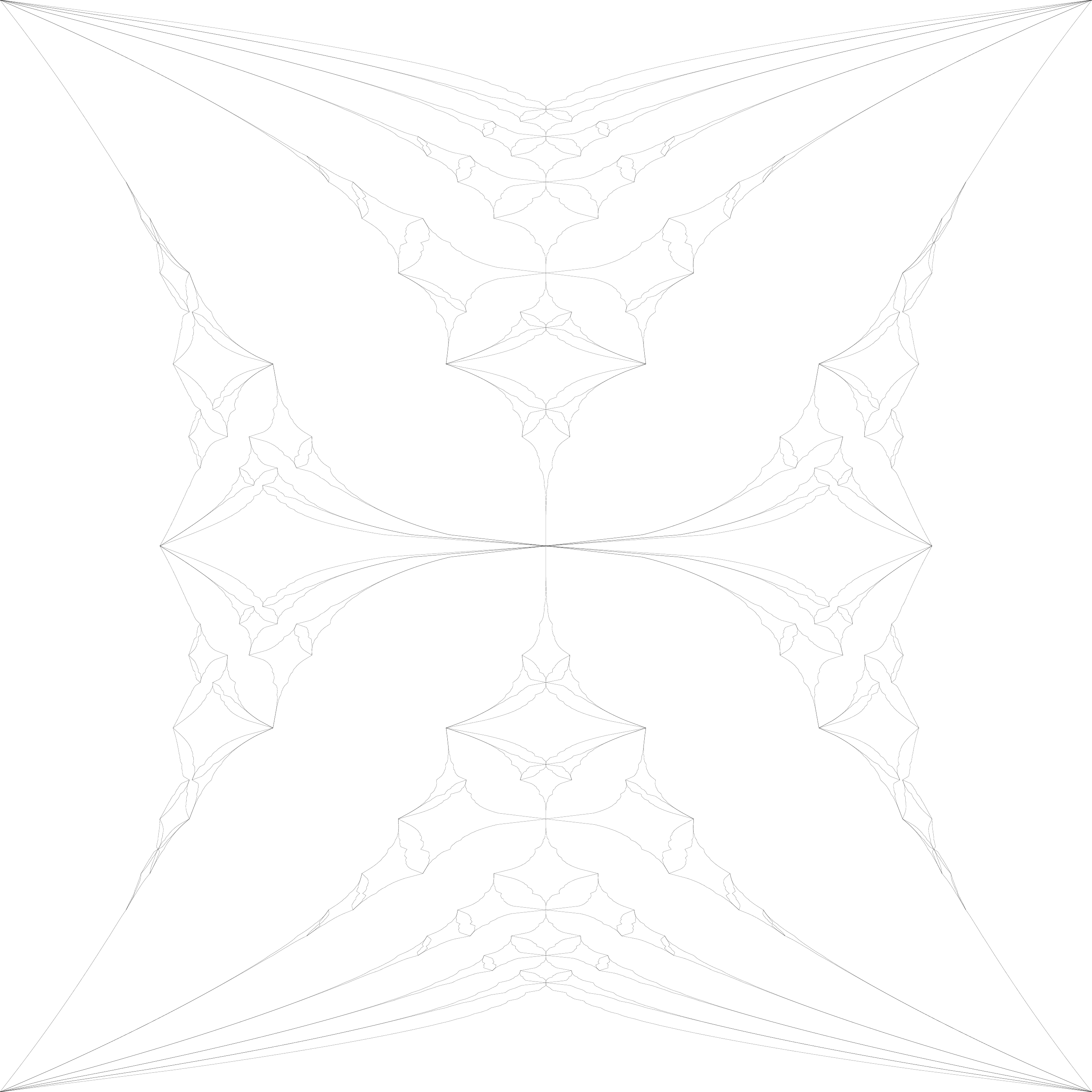} 
   \includegraphics[width=2.25in]{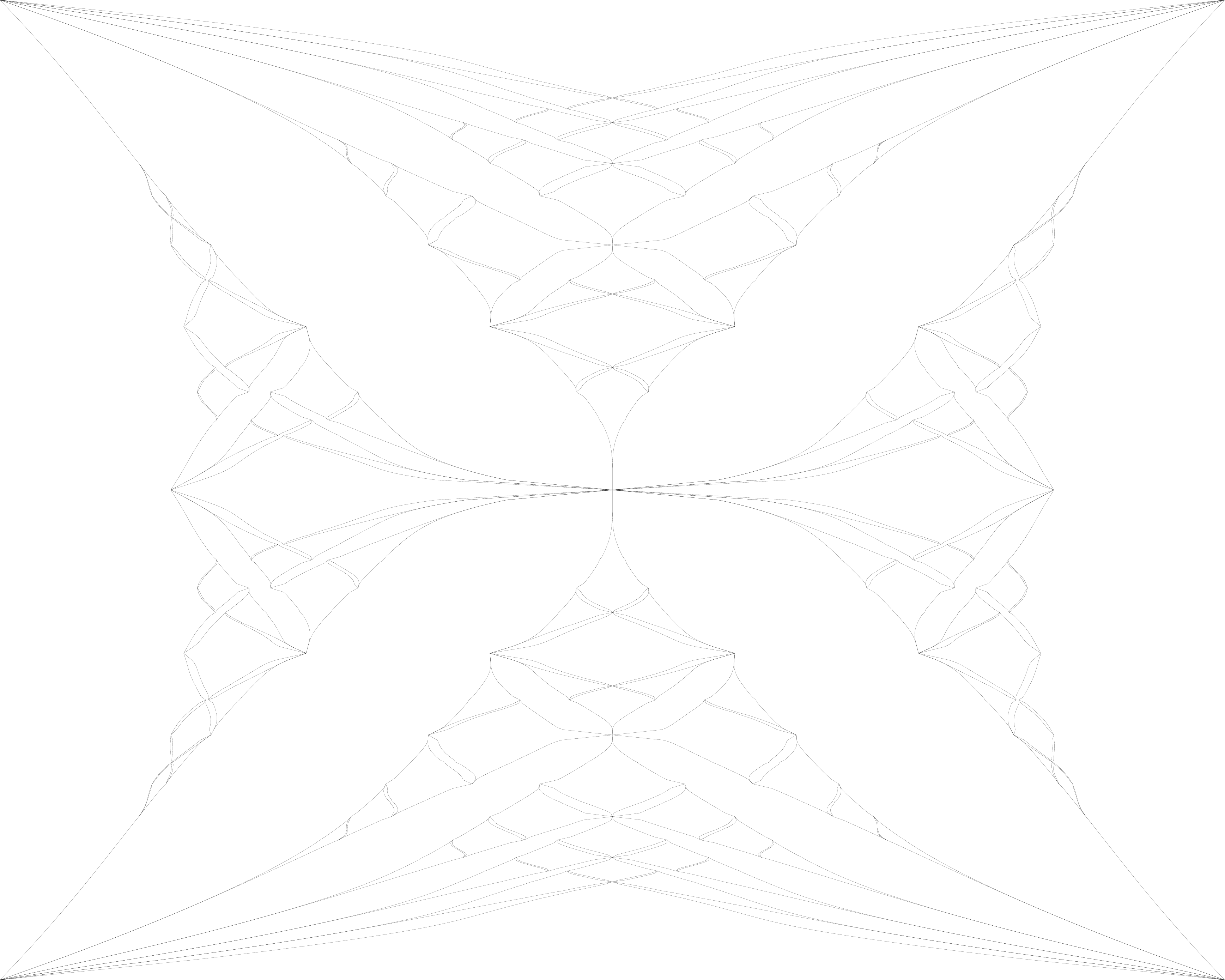} 
   \caption{Three versions with coupling constant $\lambda = 1.6, 2.0, 3.0$, respectively.}
   \label{fig:coupling}
\end{figure}

\section{Future work}

It would be useful to have a function that can select, then draw, individual pairs of butterfly wings. As things stand now, we only get half-wings of certain slopes, and it is not clear how to label them and put them all together. 

This may be related to the fractal structure of the Hofstadter butterfly, which we supsect is somehow indexed by $SL_2(\Z)$. 

\section{Summary}

We have presented a new rendering of the Hofstadter butterfly based on explicit drawings of the spectral gaps using the gap labelling method of K-theory. Discontinuities in the gap positions are revealed by this method, and we present a conjecture for identifying where the discontinuities occur, and how to ``close'' the the wing tips at the discontinuities. We have presented numerical evidence to suggest the conjecture is correct, and results in an accurate rendering of the fract-like structure for the Hofstadter butterfly.

As a side note, we observe the utility of programming our mathematical illustrations directly in the PostScript language, for optimum resolution.

\section{Acknowledgements}
This work is supported in part by an NSERC Discovery grant.

\nocite{avila06}
\nocite{bel82}
\nocite{bel90}
\nocite{brown64}
\nocite{cass05}
\nocite{choi90}
\nocite{connes94}
\nocite{gold09}
\nocite{hof76}  
\nocite{kam03} 
\nocite{lam07} 
\nocite{lam97}
\nocite{last94}
\nocite{matlab}
\nocite{puig03}
\nocite{ypma07}

\bibliography{DrawingButterflies}

\section*{Appendix 1: Code}

We list the MATLAB code that draws the butterfly wing.

There is a call to a function $Heigs(p,q,2)$ which produces the list of eigenvalues for the spectral lines of the almost Mathieu operator. This function is discussed elsewhere.

\vskip 5mm

{\footnotesize
\begin{verbatim}
% script diophant5.m
% MATLAB code to draw a range of butterfly wings
%
% We will draw a whole range of butterfly wings, from tmin to tmax

% Follow the diophantine parameterization for the gap count:
% r = t*p - s*q.

% We run through a range of t values, 1 <= t <= 5, say. 
% And 0 <= s <  t

tmin = 1; tmax = 5;
q_max = 50;  % maximum denominator we want to consider

clf; hold on  % clear the figure and set hold on

for t=tmin:tmax
    for s=0:(t-1)
        
        % first we compute the jump points
        start_pt = s/t;
        end_pt = (s+1)/t;
        jmp_pts = [start_pt, end_pt];
        for t1=1:(t-1)
            for s1=1:t1
                new_pt = (s1+s)/(t1+t);
                if (new_pt > start_pt && new_pt < end_pt)
                    if (sum(new_pt == jmp_pts)==0)  % not in list yet
                        jmp_pts = [jmp_pts new_pt];
                    end
                end
            end
        end
        
        jmp_pts = sort(jmp_pts);

        % Now we do the butterfly wings, from lower left to upper right

        for k = 1:(length(jmp_pts)-1)
            start_pt = jmp_pts(k); % set up the start and end points, plot between
            end_pt =jmp_pts(k+1);

            theta_list = []; % the list of various points to plot
            left_list = []; % the left side of the butterfly gaps
            right_list = []; % the right side of the butterfly gaps

        for q=1:q_max
            for p=0:q   % we use the fact that gcd(0,1) = gcd(1,1) = 1 to get ends
                r = t*p - s*q ;
                if (gcd(p,q)==1 && start_pt <= p/q && p/q <= end_pt )
                     if (start_pt < p/q && p/q < end_pt)
                        x = Heigs(p,q,2);
                        theta_list = [theta_list,p/q]; % add in a new theta value
                        left_list = [left_list,x(2*r)];
                        right_list = [right_list,x(2*r+1)];
                    end
                    if (p/q==start_pt)
                        x = Heigs(p,q,2);
                        theta_list = [theta_list,p/q]; % add in a new theta value
                        left_list = [left_list,x(2*r+1)];
                        right_list = [right_list,x(2*r+1)];
                    end
                    if (p/q==end_pt)
                        x = Heigs(p,q,2);
                        theta_list = [theta_list,p/q]; % add in a new theta value
                        left_list = [left_list,x(2*r)];
                        right_list = [right_list,x(2*r)];
                    end              
                end
            end
        end
        [theta_list, IX] = sort(theta_list);
        left_list = left_list(IX);
        right_list = right_list(IX);
        plot(left_list,theta_list,'-')
        plot(right_list,theta_list,'-')
        plot(-left_list,theta_list,'-')
        plot(-right_list,theta_list,'-')
        xlabel(['For r = t*p - s*q, ',num2str(tmin),'<= t <= ', num2str(tmax)])
        xlim([-4,4]);ylim([0,1]);

        end  % we end the k loop
    end % we end the s loop
end  % we end the t loop
\end{verbatim}
} 

\section*{Appendix 2: EPS file}

Here is a sample of the Encapsulated Postscript file (EPS) used to create the butterflies in the high resolution Postscript files. There is a header including the format information, then a series of Postscript commands that set up the scaling, line width, and finally a series of moveto/lineto command that draws the lines.

We use MATLAB to create the numbers, then save the list of commands directly to a file that can be used by any Postscript reader.

\vskip 5mm

{\footnotesize 
\begin{verbatim}
%!PS-Adobe-3.0 EPSF-3.0 
%%Creator: Michael P. Lamoureux, Copyright 2010  
%%Title: Almost Mathieu Butterflies 
%%CreationDate: 16-Feb-2010
%%DocumentData: Clean7Bit 
%%Origin: 0 0 
%%BoundingBox: -400 0 400 800 
100 100 scale 
.0005 setlinewidth 
1 setlinecap 
newpath 
-4.000000000000e+00 0.000000000000e+00 moveto 
-3.876300213013e+00 1.600000000000e-01 lineto 
-3.873816420419e+00 1.632653061224e-01 lineto 
...
3.625495844745e+00 7.836734693878e+00 lineto 
3.632788171307e+00 7.840000000000e+00 lineto 
4.000000000000e+00 8.000000000000e+00 lineto 
stroke 
newpath 
4.000000000000e+00 0.000000000000e+00 moveto 
3.876300213013e+00 1.600000000000e-01 lineto 
3.873816420419e+00 1.632653061224e-01 lineto 
...
-3.625495844745e+00 7.836734693878e+00 lineto 
-3.632788171307e+00 7.840000000000e+00 lineto 
-4.000000000000e+00 8.000000000000e+00 lineto 
stroke 

\end{verbatim}
} 

\section*{Appendix 3: more plots}

Now that we have the technology, let us plot two more images for looks.

The first image, Figure~17 shows 5 levels of butterflies, for a clean, simple image.
The second image, Figure~18, shows 20 levels of butterflies, for a denser, more complex image.

In all cases, the algorithm for predicting the gap discontinuities successfully removes any extraneous lines. 
\begin{figure}[hb] 
   \centering
   \includegraphics[width=6in]{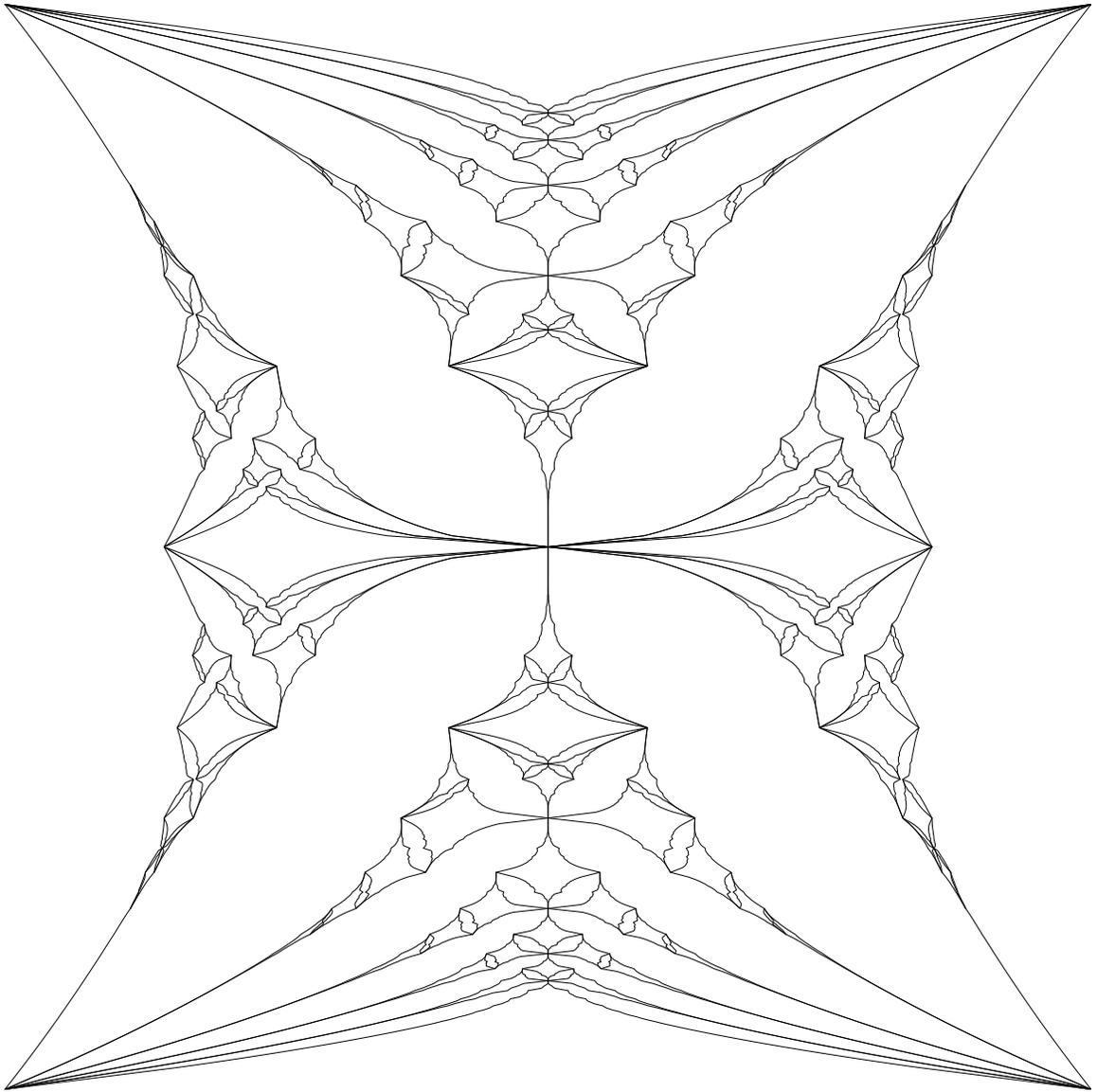} 
   \caption{Five levels of butterflies, with denominator $q$ in the range 1...50.}
   \label{fig:5levels50}
\end{figure}

.. \hfill ..
\pagebreak

\begin{figure}[h] 
   \centering
    \includegraphics[width=6in]{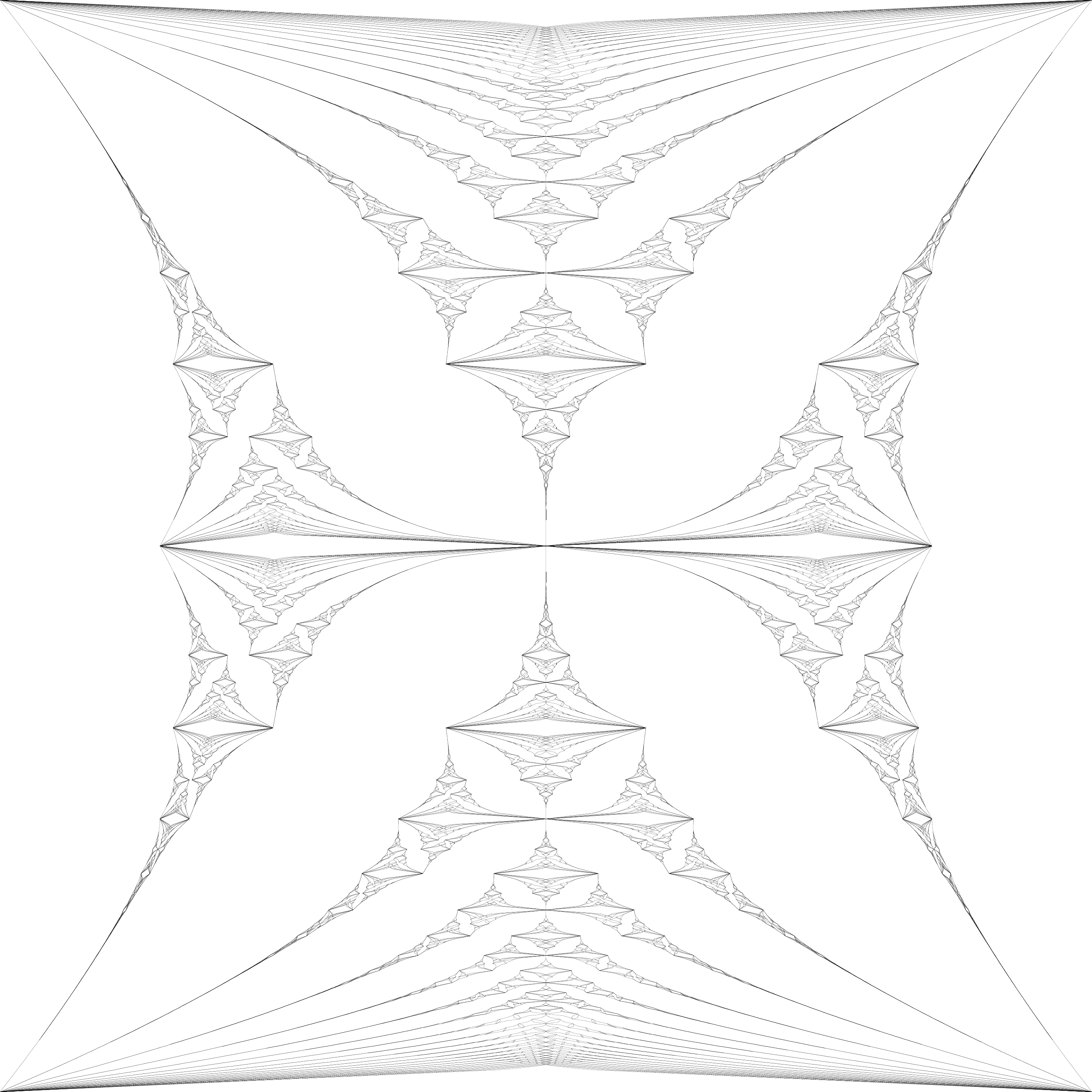} 
   \caption{Twenty levels of butterflies, with denominator $q$ in the range 1...50.}
   \label{fig:20levels100}
\end{figure}

.. \hfill ..
\pagebreak

\end{document}